\documentclass[12pt]{amsart}
\usepackage[foot]{amsaddr}
\usepackage[utf8]{inputenc}
\usepackage{mathrsfs, microtype}
\usepackage{fullpage}
\usepackage{amssymb}
\usepackage{amsmath}
\usepackage{epsfig}
\usepackage{algorithm}
\usepackage{algorithmic}
\usepackage{xfp}
\usepackage{latexsym}
\usepackage{enumitem}
\usepackage{amsmath}\usepackage{lmodern}
\usepackage{adjustbox}
\usepackage{amsfonts}
\usepackage{graphicx}
\usepackage[normalem]{ulem}
\usepackage{array}
\usepackage{colortbl}
\usepackage{bm}
\usepackage{color}\usepackage[dvipsnames]{xcolor}
\usepackage{tikz}
\usetikzlibrary{patterns, calc}
\usetikzlibrary{positioning,shapes, calc,arrows.meta}
\usepackage{xparse}

\newcommand{\jose}[1]{\mbox{}{\sf\color{red}[Jose: #1]}\marginpar{\color{blue}\Large$*$}}

\usetikzlibrary{shapes.geometric}
\usepackage{array,multirow}
\usepackage[
  bookmarks=false, 
  colorlinks,
  citecolor=brown!70!black,
  linkcolor=brown!80!black,
  urlcolor=blue!70!black,
]{hyperref}
\usepackage{relsize,exscale}

\newtheorem{cor}{Corollary}
\newtheorem{thm}{Theorem}

\newtheorem{prop}{Proposition}

\newcommand{\comm}[1]{}
 	\definecolor{lightlightgray}{rgb}{0.93, 0.93, 0.93}
 		\definecolor{llightgray}{rgb}{0.87, 0.87, 0.87}

\newcolumntype{C}[1]{>{\centering\arraybackslash }b{#1}}

\newcommand{\Sr}{{\mathcal S}}

 \def\R{17.5cm}

\usepackage{xfp}

\def\sizePoint{2.6pt}

\newcommand{\point}[2]{\fill (canvas cs:x=#1,y=#2) circle (\sizePoint); }

\newcommand\oeis[1]{\href{https://oeis.org/#1}{#1}}

\newcommand{\jluc}[1]{\mbox{}{\sf\color{blue}[Jean-Luc: #1]}\marginpar{\color{blue}\Large$*$}}

\providecommand{\subjclass}[1]{\textit{MSC:} #1}

\def\R{1} 

\newcommand{\CalcArcPoints}[2]{%
  \pgfmathsetmacro{\startX}{cos(#1)}
  \pgfmathsetmacro{\startY}{sin(#1)}
  \pgfmathsetmacro{\endX}{cos(#2)}
  \pgfmathsetmacro{\endY}{sin(#2)}
}

\newcommand{\MinMaxX}[4]{
  \ifdim #1 pt < #3 pt
    \def\minX{#1} \def\minY{#2}
    \def\maxX{#3} \def\maxY{#4}
  \else
    \def\minX{#3} \def\minY{#4}
    \def\maxX{#1} \def\maxY{#2}
  \fi
}

\newcommand{\accX}{0}
\newcommand{\accY}{0}

\newcommand{\ChainArc}[3]{%
  \CalcArcPoints{#1}{#2}%
  \MinMaxX{\startX}{\startY}{\endX}{\endY}%
  
  \pgfmathsetmacro{\dx}{\accX - \minX}
  \pgfmathsetmacro{\dy}{\accY - \minY}
  
  \begin{scope}[shift={(\dx,\dy)}]
    \draw[#3, very thick] ({cos(#1)},{sin(#1)}) arc[start angle=#1, end angle=#2, radius=\R];
  \end{scope}
  
  \pgfmathsetmacro{\accX}{\maxX + \dx}
  \pgfmathsetmacro{\accY}{\maxY + \dy}
   \pgfmathsetmacro{\accsX}{\maxX - \dx}
}

\def\NE{\ChainArc{120}{180}{green!60!black}}
\def\NEB{\ChainArc{300}{360}{brown}}
\def\SEB{\ChainArc{180}{240}{orange}}
\def\SE{\ChainArc{0}{60}{red}}
\def\E{\ChainArc{60}{120}{blue}}
\def\EB{\ChainArc{240}{300}{purple}}
\def\po{\filldraw[black] (\accX,\accY) circle (0.05);}
\newcommand{\hexsize}{1} 

\newcommand{\drawhexagon}[2]{
    \begin{scope}[shift={(#1,#2)}]
        \path[draw=black, fill=green!70!yellow]
            \foreach \i [count=\xi from 1] in {0,60,...,300} {
                ({\hexsize*cos(\i)}, {\hexsize*sin(\i)})
                \ifnum\xi<6 -- \fi
            } -- cycle;
    \end{scope}
}

\newcommand{\s}{{\texttt{width}}}
\newcommand{\ar}{{\texttt{area}}}

\newcommand{\n}{{\texttt{nbstep}}}
\newcommand{\h}{{\texttt{height}}}
\newcommand{\ki}{{\texttt{kiss}}}
\newcommand{\nbp}{{\texttt{peak}}}
\begin{document}

\title{Enumeration of Paths in a Hexagonal Circle Packing}

\author{Jean-Luc Baril$^1$}
\address{$^1$\rm LIB, Universit\'e Bourgogne Europe,
  B.P. 47 870, 21078 Dijon Cedex France}
\email{barjl@u-bourgogne.fr}

\author{Jos\'e L. Ram\'{\i}rez$^2$}
\address{$^2$Departamento de Matem\'aticas,  Universidad Nacional de Colombia,  Bogot\'a, Colombia}
\email{jlramirezr@unal.edu.co}

\date{}

\maketitle

\vspace{1em}
\noindent\textbf{Abstract.}
\quad We investigate paths in the hexagonal circle packing and enumerate them with respect to width, height, number of steps, area, and kissing number. Functional equations and the kernel method yield closed bivariate generating functions together with coefficient formulas and asymptotics. We establish bijections with skew Dyck paths, constrained Motzkin paths, and peakless Motzkin paths, and show that several of the associated counting arrays are Riordan arrays. Continued-fraction expansions for the area and kissing-number enumerators are also obtained.

\vspace{0.5cm}

\section{Introduction and motivation}

\begin{figure}[H]
\begin{tikzpicture}
    \def\r{0.5} 
    \def\n{6}   
    \def\m{5}   

    \foreach \j in {0,...,4} {
        \foreach \i in {0,...,5} {
            \pgfmathsetmacro\x{(\i + 0.5 * mod(\j,2)) * 2*\r}
            \pgfmathsetmacro\y{\j * \r * sqrt(3)}
            \draw[fill=green!15, line width=1] (\x,\y) circle(\r);
        }
    }
    \draw[dashed] (0.5,0)--(6,0);
    \draw[dashed] (0.5,0)--(0.5,4);
    \draw[fill=white,opacity=0.85,dashed] (0.5,0)--(0.5,4)-- (-0.5,4)--(-0.5,-0.5)--(0.5,-0.5)--(0.5,-0.5)--(6,-0.5)--(6,0);
    \fill (0.75,0.42) circle (2pt);
     \fill (0.5,0) circle (2pt);
\end{tikzpicture}
\begin{tikzpicture}
  \def\size{0.5} 
  \foreach \row in {0,...,4} {
    \foreach \col in {0,...,5} {
      \pgfmathsetmacro\x{2*\size*\col}
      \pgfmathsetmacro\y{\size*sqrt(3)*\row}
      \ifodd\row
        \pgfmathsetmacro\x{\x+0.14+ 0.75*\size}
      \fi
      \draw[thick] (\x,\y) 
        \foreach \i in {0,60,...,300} {
          -- ++(\i:\size)
        } -- cycle;
    }
  }
  \draw[fill=white,opacity=0.85,dashed] (0.75,0.4)--(0.75,4.4)-- (-0.25,4.4)--(-0.25,-0)--(0.5,-0)--(0.5,-0)--(6,-0)--(6,0.4);
    \draw[dashed] (0.75,0.4)--(6,0.4);   \draw[dashed] (0.75,0.4)--(0.75,4.5);
\end{tikzpicture}
\caption{The hexagonal circle packing and the associated hexagon--triangle tiling.}\label{figInt}
\end{figure}
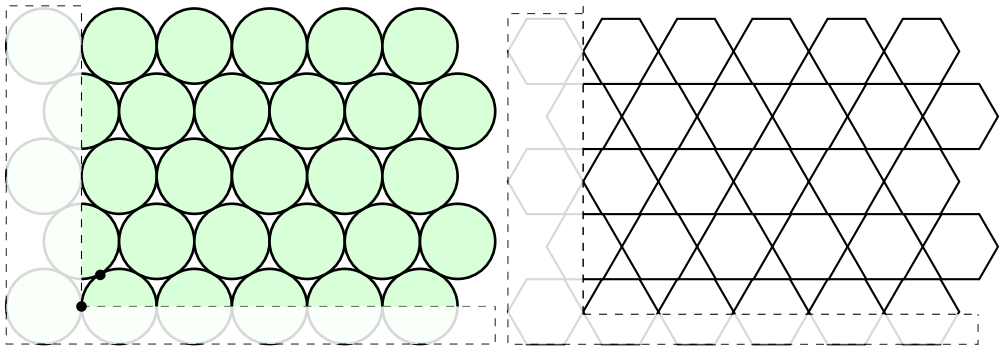

A \emph{hexagonal circle packing} is a configuration  of congruent circles in the plane arranged so that any two circles are either disjoint or tangent. Every circle is tangent to exactly six others, corresponding to its six nearest neighbors, and the set of centers of the circles forms a triangular lattice in the plane. See the left part of Figure~\ref{figInt}. The hexagonal circle packing can also be represented as a planar tiling.  A circle is surrounded by a regular hexagon determined by its  tangency points, and the hexagonal tiling  may be subdivided into equilateral triangles, see the Figure \ref{figInt}. Thus, the geometry of the packing can be studied either in terms of tangent circles or in terms of a hexagon--triangle tiling.

In this geometric setting, counting paths provides a clear mathematical way to study connectivity. It also relates to basic ideas from percolation, where one asks when long or connected paths can appear (see \cite{bousquet,essam,stauffer}).

In this paper, we study paths defined on the hexagonal circle packing. Such paths follow the arcs determined by the tangency relations between neighboring circles. 
There are three possible directions (East, North-East and South-East), and for each of them, there is a convex arc and a concave one. So, there are six possible types of arcs within the packing. These arcs are illustrated below:
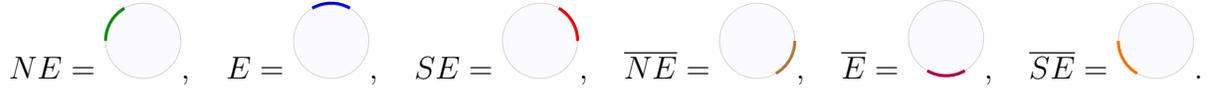
\[
 NE=\begin{tikzpicture}[scale=0.5]
        \draw[green!60!black, very thick] (120:1) arc (120:180:1); \draw[fill=blue!10,opacity=0.15 ] (0,0) circle(1);
  \end{tikzpicture}, \quad
  E=\begin{tikzpicture}[scale=0.5]
    \draw[blue, very thick] (60:1) arc (60:120:1);\draw[fill=blue!10,opacity=0.15 ] (0,0) circle(1);
  \end{tikzpicture}, \quad
   SE=\begin{tikzpicture}[scale=0.5]
    \draw[red, very thick] (0:1) arc (0:60:1);\draw[fill=blue!10,opacity=0.15 ] (0,0) circle(1);
  \end{tikzpicture}, \quad
  \overline{NE}=\begin{tikzpicture}[scale=0.5]
   \draw[brown, very thick] (300:1) arc (300:360:1);\draw[fill=blue!10,opacity=0.15 ] (0,0) circle(1);
  \end{tikzpicture}, \quad
  \overline{E}=\begin{tikzpicture}[scale=0.5]
   \draw[purple, very thick] (240:1) arc (240:300:1);\draw[fill=blue!10,opacity=0.15 ] (0,0) circle(1);
  \end{tikzpicture}, \quad
  \overline{SE}=\begin{tikzpicture}[scale=0.5]
    \draw[orange, very thick] (180:1) arc (180:240:1);\draw[fill=blue!10,opacity=0.15 ] (0,0) circle(1);
  \end{tikzpicture}.\]

For short, we denote them by $U$, $F$, $D$, $\overline{U}$, $\overline{F}$, and $\overline{D}$, respectively.  A \emph{path} in the hexagonal circle packing is formed by concatenating arcs of these six types remaining within the quarter plane $\mathbb{N}^2$. For convenience, we fix the origin at the tangency point of two horizontally aligned circles, and we set the radius of each circle to 2.

A \emph{packing path} is a path in the hexagonal circle packing  starting at the origin, ending on the $x$-axis, and consisting of steps  from  $S=\{NE, E, SE, \overline{NE}, \overline{E}, \overline{SE}\}$ (or equivalently $S=\{U, F, D, \overline{U}, \overline{F}, \overline{D}\}$). A \emph{partial packing path} is  any prefix of a packing path. See Figure~\ref{fig:path} for an illustration of a partial packing path (on the left) and a packing path (on the right).

Let $\mathcal{P}$ (resp. $\mathcal{P}'$) denote the set of all partial packing paths (resp. all packing paths). A {\it statistic} on the set $\mathcal{P}$  is a function
$\texttt{w}$ from $\mathcal{P}$ to $\mathbb{N}$. We now introduce the statistics relevant to our study. 

The {\it width} of a path $P\in\mathcal{P}$, denoted $\s(P)$, is the abscissa of its endpoint. In this setting the \emph{abscissa} of a point is defined by the horizontal displacement accumulated along the path.  Each diagonal arc of type $NE$, $SE$, $\overline{NE}$, or $\overline{SE}$ (or equivalently $U$, $D$, $\overline{U}$, $\overline{D}$) increases the abscissa by one unit, 
while each horizontal arc of type $E$ or $\overline{E}$ increases it by two units. We denote by~$\varepsilon$ the empty path, that is, the path of width~0. For $n\geq 0$, $\mathcal{P}_n$ and $\mathcal{P}_n'$ will denote the sets of paths of width $n$ in $\mathcal{P}$ and $\mathcal{P}'$, respectively. 

The {\it final height} (or simply {\it height})  of a path $P$, denoted $\h(P)$, is the unique integer such that the ordinate of the endpoint of $P$ equals $\h(P)\cdot \sqrt{3}$. Each diagonal arc $U$, $D$, $\overline{U}$, $\overline{D}$ increases the height by one unit, while each horizontal arc of type $F$ or $\overline{F}$ does not change the height. The {\it number of steps} of a path $P$, denoted $\n(P)$, is the number of arcs in the path. In the case where the path $P$ ends on the $x$-axis, we defined its {\it area}, denoted $\ar(P)$, as the number of circles lying below the path and above the line $y=-2$. The {\it kissing number} of a path $P$, denoted $\ki(P)$, is the number of circles below the path that touch it.   

For example, Figure~\ref{fig:path} (left) shows a partial packing path with $\s(P)=30$, $\h(P)=2$, and  $\n(P)=21$. Figure~\ref{fig:path} (right)  shows a  packing path with $\s(P)=28$, $\n(P)=20$, $\ar(P)=9$, and $\ki(P)=8$.


\begin{figure}[ht!]
    \centering
    \begin{tikzpicture}[scale=0.4, line cap=round]
      \draw[dashed] (0,0)--(16,0);
      \draw[dashed] (0,0)--(0,3);
      \renewcommand{\accX}{0}
      \renewcommand{\accY}{0}
      
      \filldraw[black] (0,0) circle (0.06);
      
      \NE   \filldraw[black] (\accX,\accY) circle (0.05);
      \NEB  \filldraw[black] (\accX,\accY) circle (0.05);
      \SEB  \filldraw[black] (\accX,\accY) circle (0.05);
      \EB   \filldraw[black] (\accX,\accY) circle (0.05);
      \E    \filldraw[black] (\accX,\accY) circle (0.05);
      \SE   \filldraw[black] (\accX,\accY) circle (0.05);
      \NE   \filldraw[black] (\accX,\accY) circle (0.05);
      \E    \filldraw[black] (\accX,\accY) circle (0.05);
      \EB   \filldraw[black] (\accX,\accY) circle (0.05);
      \NEB  \filldraw[black] (\accX,\accY) circle (0.05);
      \SEB  \filldraw[black] (\accX,\accY) circle (0.05);
      \EB   \filldraw[black] (\accX,\accY) circle (0.05);
      \NEB  \filldraw[black] (\accX,\accY) circle (0.05);
      \NE   \filldraw[black] (\accX,\accY) circle (0.05);
      \E    \filldraw[black] (\accX,\accY) circle (0.05);
       \EB   \filldraw[black] (\accX,\accY) circle (0.05);
       \E    \filldraw[black] (\accX,\accY) circle (0.05);
      \SE   \filldraw[black] (\accX,\accY) circle (0.05);
      \SEB  \filldraw[black] (\accX,\accY) circle (0.05);
      \EB   \filldraw[black] (\accX,\accY) circle (0.05);
      \NEB  \filldraw[black] (\accX,\accY) circle (0.05);
      \draw[fill=blue!10,opacity=0.2 ] (1,0) circle(1);
      \draw[fill=blue!10,opacity=0.2 ] (2,1.75) circle(1);
      \draw[fill=blue!10,opacity=0.2 ] (3,0) circle(1);
          \draw[fill=blue!10,opacity=0.2 ] (4,1.75) circle(1);
      \draw[fill=blue!10,opacity=0.2 ] (5,0) circle(1);
          \draw[fill=blue!10,opacity=0.2 ] (6,1.75) circle(1);
      \draw[fill=blue!10,opacity=0.2 ] (7,0) circle(1);
          \draw[fill=blue!10,opacity=0.2 ] (8,1.75) circle(1);
      \draw[fill=blue!10,opacity=0.2 ] (9,0) circle(1);
      \draw[fill=blue!10,opacity=0.2 ] (10,1.75) circle(1);
      \draw[fill=blue!10,opacity=0.2 ] (11,0) circle(1);
      \draw[fill=blue!10,opacity=0.2 ] (12,1.75) circle(1);
      \draw[fill=blue!10,opacity=0.2 ] (13,0) circle(1);
      \draw[fill=blue!10,opacity=0.2 ] (14,1.75) circle(1);
      \draw[fill=blue!10,opacity=0.2 ] (15,0) circle(1);
       \draw[fill=blue!10,opacity=0.2 ] (9,3.5) circle(1);
        \draw[fill=blue!10,opacity=0.2 ] (11,3.5) circle(1);
         \draw[fill=blue!10,opacity=0.2 ] (13,3.5) circle(1);
    \end{tikzpicture}\qquad
    \begin{tikzpicture}[scale=0.4, line cap=round]
      \draw[dashed] (0,0)--(16,0);
      \draw[dashed] (0,0)--(0,3);
      \renewcommand{\accX}{0}
      \renewcommand{\accY}{0}
      
      \filldraw[black] (0,0) circle (0.06);
      
      \NE   \filldraw[black] (\accX,\accY) circle (0.05);
      \NEB  \filldraw[black] (\accX,\accY) circle (0.05);
      \SEB  \filldraw[black] (\accX,\accY) circle (0.05);
      \EB   \filldraw[black] (\accX,\accY) circle (0.05);
      \E    \filldraw[black] (\accX,\accY) circle (0.05);
      \SE   \filldraw[black] (\accX,\accY) circle (0.05);
      \NE   \filldraw[black] (\accX,\accY) circle (0.05);
      \E    \filldraw[black] (\accX,\accY) circle (0.05);
      \EB   \filldraw[black] (\accX,\accY) circle (0.05);
      \NEB  \filldraw[black] (\accX,\accY) circle (0.05);
      \SEB  \filldraw[black] (\accX,\accY) circle (0.05);
      \EB   \filldraw[black] (\accX,\accY) circle (0.05);
      \NEB  \filldraw[black] (\accX,\accY) circle (0.05);
      \NE   \filldraw[black] (\accX,\accY) circle (0.05);
      \E    \filldraw[black] (\accX,\accY) circle (0.05);
       \EB   \filldraw[black] (\accX,\accY) circle (0.05);
       \E    \filldraw[black] (\accX,\accY) circle (0.05);
      \SE   \filldraw[black] (\accX,\accY) circle (0.05);
      \SEB  \filldraw[black] (\accX,\accY) circle (0.05);
      \SE   \filldraw[black] (\accX,\accY) circle (0.05);
      \draw[fill=blue!10,opacity=0.2 ] (1,0) circle(1);
      \draw[fill=blue!10,opacity=0.2 ] (2,1.75) circle(1);
      \draw[fill=blue!10,opacity=0.2 ] (3,0) circle(1);
          \draw[fill=blue!10,opacity=0.2 ] (4,1.75) circle(1);
      \draw[fill=blue!10,opacity=0.2 ] (5,0) circle(1);
          \draw[fill=blue!10,opacity=0.2 ] (6,1.75) circle(1);
      \draw[fill=blue!10,opacity=0.2 ] (7,0) circle(1);
          \draw[fill=blue!10,opacity=0.2 ] (8,1.75) circle(1);
      \draw[fill=blue!10,opacity=0.2 ] (9,0) circle(1);
      \draw[fill=blue!10,opacity=0.2 ] (10,1.75) circle(1);
      \draw[fill=blue!10,opacity=0.2 ] (11,0) circle(1);
      \draw[fill=blue!10,opacity=0.2 ] (12,1.75) circle(1);
      \draw[fill=blue!10,opacity=0.2 ] (13,0) circle(1);
      \draw[fill=blue!10,opacity=0.2 ] (14,1.75) circle(1);
      \draw[fill=blue!10,opacity=0.2 ] (15,0) circle(1);
       \draw[fill=blue!10,opacity=0.2 ] (9,3.5) circle(1);
        \draw[fill=blue!10,opacity=0.2 ] (11,3.5) circle(1);
         \draw[fill=blue!10,opacity=0.2 ] (13,3.5) circle(1);
    \end{tikzpicture}
    \caption{On the left, the partial packing path $U\overline{U}\overline{D}\overline{F}FDUF\overline{F}\overline{U}\overline{D}\overline{F}\overline{U}UF\overline{F}FD\overline{D}\overline{F}\overline{U}$ with 21 steps, ending at height~2 and abscissa~30. On the right, the  packing path  $U\overline{U}\overline{D}\overline{F}FDUF\overline{F}\overline{U}\overline{D}\overline{F}\overline{U}UF\overline{F}FD\overline{D}D$ with 20 steps, ending on the $x$-axis at abscissa~28, area~ 9 and kissing number~8.}
    \label{fig:path}
\end{figure}

Let us assume that the statistic $\texttt{w}$ returns either the width, or the height, or the number of steps, or the area, or the kissing number of a path. Similar statistics have been studied on the hexagonal lattice in \cite{Baram}.

For $k\geq 0$, we consider the generating function $f_k=f_k(z)$ (resp. $g_k=g_k(z)$, resp. $h_k=h_k(z)$), where the coefficient of~$z^n$ in the series expansion is the number of  paths $P\in\mathcal{P}$ such that $\texttt{w}(P)=n$, ending at height $k$ (i.e. at ordinate $k.\sqrt{3}$) with a step $U$ or $\overline{U}$ (resp., with step $D$ or $\overline{D}$, resp., with a step $F$ or $\overline{F}$).  

Let $f_k^0$ (resp. $g_k^0$, $h_k^0$) be the generating function for the paths considered in $f_k$ (resp. $g_k$, $h_k$) ending with an unbarred step.
Let $f_k^1$ (resp. $g_k^1$, $h_k^1$) be the generating function for the paths considered in $f_k$ (resp. $g_k$, $h_k$) ending with a barred step.
Obviously, we have 
$$f_k=f_k^0+f_k^1, \quad g_k=g_k^0+g_k^1, \mbox{ and } h_k=h_k^0+h_k^1, \mbox{ for any } k\geq 0.$$

 We consider the following bivariate generating functions, for $i \in \{0,1\}$:
\[
F^{i}(u):=F(x,u)=\sum_{n \geq 0} f_n^i u^n, 
\quad 
G^{i}(u):=G(x,u)=\sum_{n \geq 0} g_n^i u^n, 
\quad 
H^{i}(u):=H(x,u)=\sum_{n \geq 0} h_n^i u^n,
\]
and we set
\[
S(x,u):=F^{0}(u)+F^{1}(u)+G^{0}(u)+G^{1}(u)+H^{0}(u)+H^{1}(u).
\]
\noindent{\bf Motivation and outline of the paper.} In this paper, we provide enumerating results for partial packing paths and packing paths.  More precisely, in Sections 2,  we give bivariate generating functions and exact values for the number of these paths with respect to the width and the final height. We also provide a bijection between packing paths of width $4n$ and skew Dyck paths with $n$ steps, where the down-steps come in two colors. We establish a connection with Riordan arrays and we determine the associated $A$- and $Z$-sequences. In Section 3, we make a similar study for partial packing paths with respect to the number of steps and the final height. We provide a bijection between packing paths with $n$ steps and Motzkin paths where down-step $d$ is always followed by at least one flat-step, and each flat-step follows either a down-step or another flat-step.
In Section 4, we study packing paths with respect to the width and the area. The bivariate generating function is then expressed as a continued fraction. In Section 5, we provide the bivariate generating function for packing paths with respect to the width and the kissing number. We then focus on such paths avoiding $\overline{UD}$ and $DU$, and we establish a one-to-one correspondence with Dyck paths that make a link between  the kissing number and the number of peaks $UD$. We also prove that packing paths with a given kissing number are enumerated by the generalized Catalan numbers that count the secondary structures of RNA molecules.

\section{Partial packing of a given width}

In this section we count paths in $\mathcal{P}$ of a given width, namely those ending at a prescribed abscissa, classified by their final height.

By convention we set $f_0^0=1$, corresponding to the empty path consisting only of the origin $(0,0)$. Since a step $U$ cannot terminate at an even height, it follows that $f_{2k}^0=0$ for all $k\geq 1$.

Consider now paths  $P$ ending with a step $U$ and satisfying $\h(P)=1$. Such a path is either $P=U$, or $P=Q~U$ where $Q$ is a path of height $0$ ending with $D$. This yields \[f_1^0=x+xg^0_0.\] 
Any path $P$ ending with a step $U$ and such that $\h(P)=2k-1\geq 3$ is either $P=Q~U$ or $P=R~U$, where $Q$ (resp. $R$) is a path of height $2k-2$ ending with $\overline{U}$ (resp. ending with $D$). Hence  \[f_{2k-1}^0=x(f_{2k-2}^1+g_{2k-2}^0).\] 

Other recurrence relations for $f_k^i$, $g_k^i$, and $h_k^i$, $k\geq 0$, $i\in\{0,1\}$, can be obtained {\it mutatis mutandis}. In this way, we obtain the following  system of equations:

\begin{equation*}\left\{\begin{array}{ll}
f_0^0=1 \mbox{ and } f_{2k}^0=0, \quad & k\geq 1,\\
f_1^0=x+xg^0_0,\mbox{ and } f_{2k-1}^0=x(f_{2k-2}^1+g_{2k-2}^0), \quad &k\geq 2,\\
f_0^1=0 \mbox{ and } f_{2k-1}^1=0, \quad & k\geq 1,\\
f_{2k}^1=x(f^0_{2k-1}+h^1_{2k-1}), \quad &k\geq 1,\\
\end{array}\right.
\end{equation*}
\begin{equation*}\left\{\begin{array}{ll}
g_{2k-1}^0=0, \quad & k\geq 1,\\
g_{2k}^0=x(g^1_{2k+1}+h^0_{2k+1}), \quad &k\geq 0,\\
g_{2k}^1=0, \quad & k\geq 0,\\
g_{2k-1}^1=x(f^1_{2k}+g^0_{2k}), \quad &k\geq 1,\\
\end{array}\right.
\quad\left\{\begin{array}{ll}h_{2k}^0=0, \quad & k\geq 0,\\
h_{2k-1}^0=x^2(f^0_{2k-1}+h^1_{2k-1}), \quad &k\geq 1,\\
h_{2k}^1=0, \quad & k\geq 0,\\
h_{2k-1}^1=x^2(g^1_{2k-1}+h^0_{2k-1}), \quad &k\geq 1.\\
\end{array}\right.
\end{equation*}

Summing the previous recursions, we have
\begin{equation*}\left\{\begin{array}{ll}
F^0(u)&=1+xu(1+F^1(u)+G^0(u)),\\
F^1(u)&=xu(F^0(u)-1+H^1(u)),\\
G^0(u)&=\frac{x}{u}(G^1(u)+H^0(u)),\\
G^1(u)&=\frac{x}{u}(F^1(u)+G^0(u)-g_0^0),\\
H^0(u)&=x^2(F^0(u)-1+H^1(u)),\\
H^1(u)&=x^2(G^1(u)+H^0(u)).
\end{array}\right.
\end{equation*}
Solving the above functional equations, we deduce
$$F^0(u)={\frac {ux^3({u}^{2}{x}^{2}+1)g_0^0+2\,{u}^{3}{x}^{5}+{u}^{4}{x}^{2}+3\,
{u}^{2}{x}^{4}-{u}^{3}x+u{x}^{3}-{u}^{2}+{x}^{2}}{
{u}^{4}{x}^{2}+3\,{u}^{2}{x}^{4}-{u}^{2}+{x}^{2}}},
$$
$$F^1(u)={\frac {{u}^{2}{x}^{2} \left( 2{x}^{2}g_0^0-{u}^{2}+{x}^{2}
 \right) }{{u}^{4}{x}^{2}+3\,{u}^{2}{x}^{4}-{u}^{2}+{x}^{2}}},$$
$$G^0(u)=-{\frac {{x}^{2} \left(({u}^{2}{x}^{2}-1)g^0_0+2\,{u}^{2}{x}^{2} \right) }{{u}^{4}{x}^{2}+3\,{u}^{2}{x}^{4}-{u}^{2}+{x}^{2}}},\quad
G^1(u)=-{\frac {ux \left( ({u}^{2}{x}^{2}+2{x}^{4}-1)g^0_0+{u}
^{2}{x}^{2}+{x}^{4} \right) }{{u}^{4}{x}^{2}+3\,{u}^{2}{x}^{
4}-{u}^{2}+{x}^{2}}},$$
$$H^0(u)={\frac {{x}^{3}u \left( 2{x}^{2}g^0_0-{u}^{2}+{x}^{2} \right) 
}{{u}^{4}{x}^{2}+3\,{u}^{2}{x}^{4}-{u}^{2}+{x}^{2}}}, \quad 
H^1(u)=-{\frac {u{x}^{3} \left( ({u}^{2}{x}^{2}-1)g^0_0+2\,{u}^{2}{x}^{2} \right) }{{u}^{4}{x}^{2}+3\,{u}^{2}{x}^{4}-{u}^{2}+{x}^{2}}}.
$$
In order to compute $g_0^0=G^0(0)$, we use the kernel method~\cite{ban, banf, pro, pro2, pro3}
on $F^1(u)$. This method consists in cancelling the denominator of $F^1(u)$
by finding $u^2$ as an algebraic function $r$ of~$x$. Therefore,  if we substitute $u^2$ with $r$ in the numerator, then it necessarily vanishes (in order to counterbalance the cancellation of the denominator), which determines the value of $g_0^0$.
Thus, we factorize the denominator ${u}^{4}{x}^{2}+3\,{u}^{2}{x}^{4}-{u}^{2}+{x}^{2}=x^2(u^2-r)(u^2-s)$ with
 $$r={\frac {1-3\,{x}^{4}-\sqrt {9\,{x}^{8}-10\,{x}^{4}+1}}{2{x}^{2}}},\quad s={\frac {1-3\,{x}^{4}+\sqrt {9\,{x}^{8}-10\,{x}^{4}+1}}{2{x}^{2}}}.$$

The cancellation of the numerator $2x^3g^0_0-r^2+x^2$ gives 
$$g^0_0=\frac{r-x^2}{2x^2}={\frac {1-5\,{x}^{4}-\sqrt {9\,{x}^{8}-10\,{x}^{4}+1}}{4{x}^{4}}}.$$

Now, substituting $g_0^0$ with this value in the  generating functions, and simplifying by $(u^2-r)$ in both the numerator and  the denominator, we obtain the following result.

\begin{thm}\label{teo1} We have 
$$F^0(u)=1+{\frac {u\,(2-3\,x^4-r\,x^2 )}{2\,x \left(s-{u}^{2} \right) }}
,
\quad  F^1(u)=\frac{u^2}{s-{u}^{2} },
\quad G^0(u)=\frac{3x^2+r}{2\left(s-{u}^{2}\right)},  \mbox{ and}$$

$$G^1(u)=\frac{u(x^2+r)}{2x\left(s-{u}^{2} \right)},\quad H^0(u)=\frac{xu}{s-{u}^{2} },
\quad H^1(u)=\frac{(3\,x^2+r)\,x\,u}{2\left(s-{u}^{2} \right)}.$$    
Finally, the bivariate generating function $S(x,u)$, where the coefficient of~$x^nu^k$ is the number of paths of width $n$ ending at height $k$, satisfies \begin{align*}S(x,u)&={\frac {u \left( {\it s}\,{x}^{2}-2\,{x}^{2}-1 \right) -{\it s}
\,{x}^{3}-x}{2{x}^{3} \left( {u}^{2}-{\it s} \right) }}
\\
&={\frac {3\,u{x}^{4}-3\,{x}^{5}+4\,{x}^{2}u+u+3\,x-(u-x)\sqrt {9\,{x}^{8}-
10\,{x}^{4}+1}}{2x \left( 1-2\,
{u}^{2}{x}^{2}-3\,{x}^{4}+\sqrt {9\,{x}^{8}-10\,{x}^{4}+1}\right) }}.
\end{align*}
\end{thm}

The first terms of the series expansion of $S(x,u)$ are
\begin{align*}    
1+xu+u^2x^2+&(u^3+2u)x^3+(u^4+2)x^4+(u^5+2u^3+4u)x^5+(u^6+6u^2)x^6+\\
&(u^7+2u^5+8u^3+10u)x^7+(u^8+10u^4+\bm{10})x^8+\\&(u^9+2u^7+12u^5+18u^3+20u)x^9+
(u^{10}+14u^6+38u^2)x^{10}+\\&(u^{11}+2u^9+16u^7+26u^5+56u^3+58u)x^{11}+O(x^{12}).
\end{align*}

\begin{figure}[ht!]
    \centering
\begin{tikzpicture}[scale=0.4, line cap=round]
\draw[dashed] (0,0)--(4,0);
\renewcommand{\accX}{0}
\renewcommand{\accY}{0}
\po\NE\po\E\po\SE\po\NE\po\E\po\SE\po
\filldraw[black] (0,0) circle (0.02);
\filldraw[black] (\accX,\accY) circle (0.02);
\draw[fill=blue!10,opacity=0.2 ] (1,0) circle(1);
\draw[fill=blue!10,opacity=0.2 ] (3,0) circle(1);
\end{tikzpicture}\quad
\begin{tikzpicture}[scale=0.4, line cap=round]
\draw[dashed] (0,0)--(4,0);
\renewcommand{\accX}{0}
\renewcommand{\accY}{0}
\po\NE\po\E\po\SE\po\NE\po\NEB\po\SEB\po\SE\po
\filldraw[black] (0,0) circle (0.02);
\filldraw[black] (\accX,\accY) circle (0.02);
\draw[fill=blue!10,opacity=0.2 ] (1,0) circle(1);
\draw[fill=blue!10,opacity=0.2 ] (3,0) circle(1);
\end{tikzpicture}\quad
\begin{tikzpicture}[scale=0.4, line cap=round]
\draw[dashed] (0,0)--(4,0);
\renewcommand{\accX}{0}
\renewcommand{\accY}{0}
\po\NE\po\E\po\EB\po\E\po\SE\po
\filldraw[black] (0,0) circle (0.02);
\filldraw[black] (\accX,\accY) circle (0.02);\draw[fill=blue!10,opacity=0.2 ] (1,0) circle(1);
\draw[fill=blue!10,opacity=0.2 ] (3,0) circle(1);
\end{tikzpicture}\quad
\begin{tikzpicture}[scale=0.4, line cap=round]
\draw[dashed] (0,0)--(4,0);
\renewcommand{\accX}{0}
\renewcommand{\accY}{0}
\po\NE\po\E\po\EB\po\NEB\po\SEB\po\SE\po
\filldraw[black] (0,0) circle (0.02);
\filldraw[black] (\accX,\accY) circle (0.02);\draw[fill=blue!10,opacity=0.2 ] (1,0) circle(1);
\draw[fill=blue!10,opacity=0.2 ] (3,0) circle(1);
\end{tikzpicture}\quad
\begin{tikzpicture}[scale=0.4, line cap=round]
\draw[dashed] (0,0)--(4,0);
\renewcommand{\accX}{0}
\renewcommand{\accY}{0}
\po\NE\po\NEB\po\SEB\po\SE\po\NE\po\E\po\SE\po
\filldraw[black] (0,0) circle (0.02);
\filldraw[black] (\accX,\accY) circle (0.02);\draw[fill=blue!10,opacity=0.2 ] (1,0) circle(1);
\draw[fill=blue!10,opacity=0.2 ] (3,0) circle(1);
\end{tikzpicture}\\[0.1cm]
\begin{tikzpicture}[scale=0.4, line cap=round]
\draw[dashed] (0,0)--(4,0);
\renewcommand{\accX}{0}
\renewcommand{\accY}{0}
\po\NE\po\NEB\po\SEB\po\SE\po\NE\po\NEB\po\SEB\po\SE\po
\filldraw[black] (0,0) circle (0.02);
\filldraw[black] (\accX,\accY) circle (0.02);\draw[fill=blue!10,opacity=0.2 ] (1,0) circle(1);
\draw[fill=blue!10,opacity=0.2 ] (3,0) circle(1);
\end{tikzpicture}\quad
\begin{tikzpicture}[scale=0.4, line cap=round]
\draw[dashed] (0,0)--(4,0);
\renewcommand{\accX}{0}
\renewcommand{\accY}{0}
\po\NE\po\NEB\po\SEB\po\EB\po\E\po\SE\po
\filldraw[black] (0,0) circle (0.02);
\filldraw[black] (\accX,\accY) circle (0.02);\draw[fill=blue!10,opacity=0.2 ] (1,0) circle(1);
\draw[fill=blue!10,opacity=0.2 ] (3,0) circle(1);
\end{tikzpicture}\quad
\begin{tikzpicture}[scale=0.4, line cap=round]
\draw[dashed] (0,0)--(4,0);
\renewcommand{\accX}{0}
\renewcommand{\accY}{0}
\po\NE\po\NEB\po\SEB\po\EB\po\NEB\po\SEB\po\SE\po
\filldraw[black] (0,0) circle (0.02);
\filldraw[black] (\accX,\accY) circle (0.02);\draw[fill=blue!10,opacity=0.2 ] (1,0) circle(1);
\draw[fill=blue!10,opacity=0.2 ] (3,0) circle(1);
\end{tikzpicture}\quad
\begin{tikzpicture}[scale=0.4, line cap=round]
\draw[dashed] (0,0)--(4,0);
\renewcommand{\accX}{0}
\renewcommand{\accY}{0}
\po\NE\po\NEB\po\NE\po\E\po\SE\po\SEB\po\SE\po
\filldraw[black] (0,0) circle (0.02);
\filldraw[black] (\accX,\accY) circle (0.02);\draw[fill=blue!10,opacity=0.2 ] (1,0) circle(1);
\draw[fill=blue!10,opacity=0.2 ] (3,0) circle(1);\draw[fill=blue!10,opacity=0.2 ] (2,1.75) circle(1);
\end{tikzpicture}\quad
\begin{tikzpicture}[scale=0.4, line cap=round]
\draw[dashed] (0,0)--(4,0);
\renewcommand{\accX}{0}
\renewcommand{\accY}{0}
\po\NE\po\NEB\po\NE\po\NEB\po\SEB\po\SE\po\SEB\po\SE\po
\filldraw[black] (0,0) circle (0.02);
\filldraw[black] (\accX,\accY) circle (0.02);\draw[fill=blue!10,opacity=0.2 ] (1,0) circle(1);
\draw[fill=blue!10,opacity=0.2 ] (3,0) circle(1);\draw[fill=blue!10,opacity=0.2 ] (2,1.75) circle(1);
\end{tikzpicture}
    \caption{The 10 paths ending on the $x$-axis at abscissa 8.}
    \label{fig:pathSec1}
\end{figure}
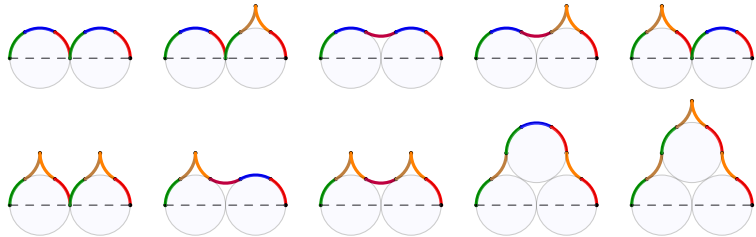

In Figure~\ref{fig:pathSec1} we show the 10  paths of width~8 ending on the $x$-axis, 
corresponding to the bold coefficient in the above series.

Using the Vieta relations $rs=1$ and $r+s=\frac{1-3x^4}{x^2}$, we deduce the following.

\begin{cor}\label{cor1} We have 
$$[u^{2k-1}]F^0(u)=\frac{2-rx^2-3x^4}{2x}r^{k}, ~k\geq 1,\quad  [u^{2k}]F^1(u)=r^{k}, ~k\geq 1,$$
$$[u^{2k}]G^0(u)=\frac{3x^2+r}{2}r^{k+1}, ~k\geq 1, \quad 
[u^{2k+1}]G^1(u)=\frac{x^2+r}{2x}r^{k+1}, ~k\geq 0, \
$$ $$[u^{2k-1}]H^0(u)=x\,r^k,~k\geq 1,\quad [u^{2k-1}]H^1(u)=\frac{3x^2+r}{2}\,r^k,~k\geq 1,$$ $$[u^{2k}]S(x,u)=\frac{sx^2+1}{2x^2}r^{k+1}
,~k\geq 0, \qquad [u^{2k+1}]S(x,u)=\frac{1+2x^2-sx^2}{2x^3}r^{k+1}
,~k\geq 0, \mbox{ and}$$    
all other coefficients vanish, except $[u^0]F^0(u)=1$.
\end{cor}

In the following theorem we give an expression to compute the  coefficient $x^nu^k$ in the series expansion of $S(x,u)$. 

\begin{thm} Let $s(n,k)$ denote the number of paths of width $n$ ending at height $k$. Then 
\begin{align*} s(4n-2k,2k)&=\frac{1}{2}\left(w(4n-2k,k)+w(4n-2(k-1),k+1)\right),\\
s(4n-2k-1,2k+1)&=\frac{1}{2}\left(w(4n-2(k-1),k+1)-w(4n-2k,k)\right), \\
s(4n-2k+1,2k+1)&=w(4n-2(k-1),k+1),
\end{align*}
and  $s(n,k)=0$ otherwise, where 
\[
w(4n-2k,k)=\frac{k}{n}\sum_{\ell=0}^{n-k} 3^{\,n-k-\ell}\binom{n+\ell-1}{\ell}\binom{n}{n-k-\ell}, \, n\geq 1, k\geq 0.
\]
    \end{thm}
 \begin{proof}
Define
\[
\hat{r}=x^{1/2}\,r(x^{1/4})
=\frac{1-3x-\sqrt{\,1-10x+9x^2\,}}{2},
\]
which  satisfies the  functional equation $\hat{r}=x\Phi(\hat{r})$, where 
$$\Phi(u)=\frac{1+3u}{1-u}.$$
By Lagrange inversion (see, e.g., \cite{merl}), for integers $n, k \geq 1$ we obtain
\begin{align*}
[x^n]\hat{r}^k&=\frac{k}{n}[x^{n-k}]\Phi(x)^n\\
&=\frac{k}{n}[x^{n-k}]\left(\frac{1+3x}{1-x}\right)^n\\
&=\frac{k}{n}\sum_{\ell=0}^{n-k}3^{n-k-\ell}\binom{n + \ell - 1}{\ell}\binom{n}{n-k-\ell}.
\end{align*}
Thus, for $n\geq 1$ we have 
$$[x^{4n-2k}]r^k=\frac{k}{n}\sum_{\ell=0}^{n-k}3^{n-k-\ell}\binom{n + \ell - 1}{\ell}\binom{n}{n-k-\ell}, $$
while $[x^m]r^k=0$ if $m\not\equiv 2k \pmod{4}$. 

Using Corollary~\ref{cor1} and setting $w(4n-2k,k):=[x^{4n-2k}]r^k$, we then obtain explicit closed forms for the numbers  $f^0(n,k)$,  $f^1(n,k)$,  $g^0(n,k)$, $g^1(n,k)$,  $h^0(n,k)$, and $h^1(n,k)$ which count partial packing paths of width $n$ ending at height $k$ whose last step is of type
$U$, $\overline{U}$, $F$, $\overline{F}$, $D$, and $\overline{D}$, respectively.
 \end{proof}

\begin{cor}The generating function $S(x)=S(x,1)$, where the coefficient of~$x^n$ is the total number of paths of width $n$, satisfies
$$S(x)=\frac{-1 - x + 4 x^3 + 9 x^4 - 3 x^5 + (1 + x)\sqrt{1 - 10 x^4 + 9 x^8}}{4 x^3 (1 - 3 x^2)}.$$
\end{cor}
The first terms of the series expansion are 
$$1+x+x^2+3x^3+3x^4+7x^5+7x^6+21x^7+21x^8+53x^9+53x^{10}+159x^{11}+O(x^{12}).$$
Moreover, for every $n\geq 1$,  the coefficients of  $x^{2n}$ and  $x^{2n-1}$ coincide, and the sequence of coefficients of $x^{2n}$ for $n\geq 0$ corresponds to \oeis{A368773}, which currently appears in the OEIS without a known combinatorial interpretation.

\begin{cor} 
The generating function whose coefficient of~$x^n$ is the number of paths of width $n$ ending on the $x$-axis satisfies
$$S(x,0)=\frac{x^2+r}{2x^2}=\frac{1-x^{4}-\sqrt{1-10x^{4}+9x^{8}}}{4 x^{4}}.$$
Moreover, we have $[x^0]S(x,0)=1$, $$[x^{4n}]S(x,0)=\frac{2}{n}\sum\limits_{k=0}^n4^k{n\choose k}{n\choose k+1}, \mbox{ for } n\geq 1,$$   and 0 otherwise.  
An asymptotic approximation of the $4n$-th term is 
$$ \frac{3^{1+2n}}{\sqrt{\pi}(1+2n)^{3/2}}.$$
\end{cor}

The first terms of the series expansion are 
$$1+2x^4+10x^8+58x^{12}+370x^{16}+2514x^{20}+17850x^{24}+O(x^{28}).$$
The sequence of coefficients of $x^{4n}$ for $n\geq 0$ corresponds to the sequence \oeis{A086871}, which counts skew Dyck paths of semilength $n$ with  down-steps in two colors.

By calculating $\partial_u(S(x,u))\rvert_{u=1}$, and using classical methods (see~\cite{flajolet,orlov}) for the asymptotic approximation of the coefficient of~$x^n$, we obtain the following result.

\begin{cor} An asymptotic for the expected height in all paths of a given width is given by $$\sqrt{\frac{\pi n}{2}}.$$
\end{cor}

\subsection{A bijection with the skew Dyck paths} 
We conclude this section by presenting a constructive bijection $f$ between packing paths of width $4n$ ending on the $x$-axis and skew Dyck paths with $n$ steps, where the down-steps come in two colors. 

Recall that a \emph{skew Dyck path} is a lattice path in the quarter plane starting at the origin, ending on the $x$-axis, and consisting of steps $u=(1,1)$, $d=(1,-1)$ and $\ell=(-1,-1)$. Let $\mathcal{K}_n$ denote the set of skew Dyck paths of length $n$ in which each down-step can be taken in one of two colors, say $d$ and $\textcolor{red}{d'}$.  We also set $\mathcal{K}=\cup_{n\geq 0}\mathcal{K}_n$.

\begin{prop} Define a map $f:\mathcal{P}'\to\mathcal{K}$ recursively as follows.  For the empty path we set $f(\epsilon)=\epsilon$.  
For $P\neq \epsilon$, we define
$$f(P)=\left\{ \begin{array}{lll}
ud~f(Q)&\mbox{ if } P=UFD~ Q \quad &\mbox{ with } Q\in \mathcal{P}',\\
u\textcolor{red}{d'} ~f(Q)&\mbox{ if } P=U\overline{U}\overline{D}D ~Q\quad &\mbox{ with } Q\in \mathcal{P}',\\
u~f(U~Q)~\ell&\mbox{ if } P=UF~\overline{F}~Q\quad &\mbox{ with } UQ\in \mathcal{P}',\\
u~f(Q)~d &\mbox{ if } P=U\overline{U}~ Q~ \overline{D}D\quad &\mbox{ with } Q\in \mathcal{P}', Q\neq \epsilon,\\
u~f(U~R)~\textcolor{red}{d'} ~f(Q)&\mbox{ if } U\overline{U}~ Q ~\overline{D}\overline{F} ~R &\mbox{ with } Q\in \mathcal{P}',~UR\in \mathcal{P}'.\\
\end{array}
\right.$$
Then $f$ is a bijection from $\mathcal{P}'$ onto $\mathcal{K}$.
\end{prop}

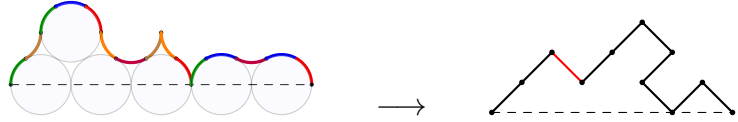
\begin{figure}[H] 
    \centering
\begin{tikzpicture}[scale=0.4, line cap=round]
\draw[dashed] (0,0)--(10,0);
\renewcommand{\accX}{0}
\renewcommand{\accY}{0}
\po\NE\po\NEB\po\NE\po
\E\po\SE\po\SEB\po
\EB\po\NEB\po\SEB\po\SE
\po\NE\po
\E\po\EB\po\E\SE\po
\filldraw[black] (0,0) circle (0.02);
\filldraw[black] (\accX,\accY) circle (0.02);
\draw[fill=blue!10,opacity=0.2 ] (1,0) circle(1);
\draw[fill=blue!10,opacity=0.2 ] (2,1.75) circle(1);
\draw[fill=blue!10,opacity=0.2 ] (3,0) circle(1);
\draw[fill=blue!10,opacity=0.2 ] (5,0) circle(1);
\draw[fill=blue!10,opacity=0.2 ] (7,0) circle(1);
\draw[fill=blue!10,opacity=0.2 ] (9,0) circle(1);
\end{tikzpicture}\qquad $\longrightarrow$\qquad
\begin{tikzpicture}[scale=0.4, line cap=round]
\draw[dashed] (0,0)--(8,0);
\draw[line width=0.8pt] (0,0)--(1,1);
\draw[line width=0.8pt] (1,1)--(2,2);
\draw[red,line width=0.8pt] (2,2)--(3,1);
\draw[line width=0.8pt] (3,1)--(4,2)--(5,3)--(6,2)--(5,1)--(6,0)--(7,1)--(8,0);
\point{0}{0}
\point{1cm}{1cm}\point{2cm}{2cm}\point{3cm}{1cm}\point{4cm}{2cm}\point{5cm}{3cm}\point{6cm}{2cm}\point{5cm}{1cm}\point{6cm}{0cm}\point{7cm}{1cm}\point{8cm}{0cm}
\end{tikzpicture}
  \caption{The packing path $P=U\overline{U}~UFD~\overline{D}D\overline{F}U\overline{D}DUF\overline{F}FD$ and its image by $f$: 
$f(P)=uu\textcolor{red}{d'}uud\ell dud$.}
    \label{bijection}
\end{figure}

Due to the recursive definition, the image under $f$ of a packing path in $\mathcal{P}'_{4n}$ is a skew Dyck path in $\mathcal{K}_n$, and it is not difficult to check the injectivity of $f$. The cardinality of $\mathcal{P}'_{4n}$ and $\mathcal{K}_n$ being the same, we deduce the bijectivity of $f$.  We refer to Figure~\ref{bijection} for an illustration of $f$ when $P=U\overline{U}~UFD~\overline{D}D\overline{F}U\overline{D}DUF\overline{F}FD$ and
$f(P)=uu\textcolor{red}{d'}uud\ell dud$.

\subsection{A connection with Riordan arrays}
Consider the matrices $\mathcal{S}_e=[s(4n-2k,2k)]_{n, k\geq 0}$ and  
$\mathcal{S}_o=[s(2n+1,2k+1)]_{n, k\geq 0}$.  The first few rows of $\mathcal{S}_e$ and $\mathcal{S}_o$ are
\[
\footnotesize 
\mathcal{S}_e=\left(
\begin{array}{ccccccc}
 1 & 0 & 0 & 0 & 0 & 0 & 0 \\
 2 & 1 & 0 & 0 & 0 & 0 & 0 \\
 10 & 6 & 1 & 0 & 0 & 0 & 0 \\
 58 & 38 & 10 & 1 & 0 & 0 & 0 \\
 370 & 254 & 82 & 14 & 1 & 0 & 0 \\
 2514 & 1774 & 642 & 142 & 18 & 1 & 0 \\
 17850 & 12822 & 4986 & 1286 & 218 & 22 & 1 \\
\end{array} 
\right) 
\quad 
\mathcal{S}_o=\left(
\begin{array}{ccccccc}
 1 & 0 & 0 & 0 & 0 & 0 & 0 \\
 2 & 1 & 0 & 0 & 0 & 0 & 0 \\
 4 & 2 & 1 & 0 & 0 & 0 & 0 \\
 10 & 8 & 2 & 1 & 0 & 0 & 0 \\
 20 & 18 & 12 & 2 & 1 & 0 & 0 \\
 58 & 56 & 26 & 16 & 2 & 1 & 0 \\
 116 & 138 & 108 & 34 & 20 & 2 & 1 \\
\end{array}
\right).
\]

Both matrices are Riordan arrays. Recall that a \emph{Riordan array} is an infinite lower triangular matrix defined by a pair of formal power series $(g(x), f(x))$ with $g(0)\neq 0$, $f(0)=0$ and $f'(0)\neq 0$, 
such that the $k$th column has generating function $g(x)f(x)^k$. If we multiply $(g, f)$ by a column vector $(c_0, c_1, \dots)^T$
with the generating function  $h(x)$, then the resulting column vector has generating function $g(x)h(x)(f(x))$.

The product of two Riordan arrays $(g(x),f(x))$ and $(h(x),\ell(x))$ is defined by \[
(g(x),f(x)) \ast (h(x),\ell(x))=\bigl(g(x)\,h(f(x)),\, \ell(f(x))\bigr).
\]
It is well known that the set of Riordan arrays forms a group under this operation \cite{Riordan2}. The identity element is $I=(1,x)$, and the inverse of $(g(x),f(x))$ is given by
\begin{align}\label{inRiordan}
(g(x), f(x))^{-1}=\left(\frac{1}{(g\circ \overline{f})(x)}, \, \overline{f}(x)\right),
\end{align}
where $\overline{f}(x)$ denotes the compositional inverse of $f(x)$.

Notice that the full matrix $\mathcal{S}=[s(n,k)]_{n, k\geq 0}$ is not a Riordan array, but rather a double Riordan array (cf. \cite{BookRiordan}). From Corollary~\ref{cor1} and the definition of a Riordan array, we obtain the following result.  

\begin{thm}\label{Riordan}
The matrices $\Sr_e$ and $\Sr_o$ are Riordan arrays. More precisely,
$$\Sr_e=\left(\frac{1 - x -\sqrt{1 - 10 x + 9 x^2}}{4x}, \frac{1 - 3x -\sqrt{1 - 10 x + 9 x^2}}{2} \right)$$
and
$$\Sr_o=\left(\frac{1 + 2 x - 5 x^2 - 6 x^3 - (1 + 2 x)\sqrt{1 - 10 x^2 + 9 x^4}}{4x^3}, \frac{1 - 3 x^2 -\sqrt{1 - 10 x^2 + 9 x^4}}{2x} \right).$$
\end{thm}

Riordan arrays admit a useful structural characterization in terms of two auxiliary sequences, known as the $A$-sequence and the $Z$-sequence. These sequences encode, respectively, the linear relations between consecutive rows  and together with the initial value $d_{0,0}$ they completely determine the array. More precisely, an infinite lower triangular matrix $\mathcal{D} = [d_{n,k}]_{n,k\in\mathbb{N}}$ is a Riordan array if and only if $d_{0,0}\neq 0$ and there exist sequences $A=(a_0\neq 0,a_1,a_2,\dots)$ and $Z=(z_0,z_1,z_2,\dots)$ such that
\begin{align*}
d_{n+1,k+1}&=a_0 d_{n,k} + a_1 d_{n,k+1} + a_2 d_{n,k+2} + \cdots, \quad n,k\geq 0, \\
d_{n+1,0}&=z_0 d_{n,0} + z_1 d_{n,1} + z_2 d_{n,2} + \cdots, \quad n\geq 0.
\end{align*}

Moreover, if $\mathcal{D}=(g(x),f(x))$ is a Riordan array with inverse 
$\mathcal{D}^{-1}=(d(x),h(x))$, then the generating functions of the $A$- and $Z$-sequences are given by
\[
A(x)=\frac{x}{h(x)}, 
\qquad 
Z(x)=\frac{1}{h(x)}\bigl(1-d_{0,0}\, d(x)\bigr).
\]
Further details on these characterizations may be found in \cite{BookRiordan}

In Corollary~\ref{AZSeq} we determine the $A$-sequence and $Z$-sequence for the Riordan array $\Sr_e$.

\begin{cor}\label{AZSeq}
For any positive integers $n, k$, we have the following relations 
\begin{align*}
s(4n-2k+2,2k+2)&=s(4 n - 2 k,2k)+4\sum_{\ell\geq 0}s(4n-2(k+\ell)-2,2k+2\ell+2), \\
s(4n-2k+2,0)&=2\sum_{\ell\geq 0}s(4n-2(k+\ell)-2,4\ell)  + 6\sum_{\ell\geq 0}s(4n-2(k+\ell)-2,4\ell+2).
\end{align*}
\end{cor}
\begin{proof}
From Theorem~\ref{Riordan} and~\eqref{inRiordan} we obtain
$$\Sr_e^{-1}=\left(\frac{1-x}{1+x}, \frac{x - x^2}{1+3x} \right).$$
Hence, the generating functions of the $A$- and $Z$-sequences of the Riordan array $\Sr_e$ are
\[A(x)=\frac{1 + 3 x}{1 - x} \quad \text{and} \quad Z(x)=\frac{2 + 6 x}{1 - x^2}.\]
Therefore, the corresponding sequences are
$$A=(1,4,4,4,\dots) \quad \text{and} \quad Z=(2,6,2,6,\dots).$$
Combining the above relations we obtain the desired result. 
\end{proof}
For the Riordan array $\Sr_o$, it is possible to obtain similar expressions; however, they are quite involved, so we omit them here for brevity.


\section{Partial packing paths with a given number of steps}

In this section, we count partial packing paths in $\mathcal{P}$ with respect to the number of steps and the final height.

Since the steps $U$, $\overline{U}$, $D$, and $\overline{D}$ each have width 1, the recursive formulas for $f_k^0$, $f_k^1$, $g_k^0$, and $g_k^1$ ($k\geq 0$) are identical to those presented in the previous section. Here, however, the steps $F$ and $\overline{F}$ are counted as a single steps, unlike in the previous section where they were assigned width 2. This modification leads to the following new relations for $h_k^0$ and $h_k^1$.

\begin{equation*}\left\{\begin{array}{ll}h_{2k}^0=0, \quad & k\geq 0,\\
h_{2k-1}^0=x(f^0_{2k-1}+h^1_{2k-1}), \quad &k\geq 1,\\
h_{2k}^1=0, \quad & k\geq 0,\\
h_{2k-1}^1=x(g^1_{2k-1}+h^0_{2k-1}), \quad &k\geq 1.\\
\end{array}\right.
\end{equation*}
Summing these recursions, we obtain the following system of functional equations:
\begin{equation*}\left\{\begin{array}{ll}
F^0(u)&=1+xu(1+F^1(u)+G^0(u)),\\
F^1(u)&=xu(F^0(u)-1+H^1(u)),\\
G^0(u)&=\frac{x}{u}(G^1(u)+H^0(u)),\\
G^1(u)&=\frac{x}{u}(F^1(u)+G^0(u)-g_0^0),\\
H^0(u)&=x(F^0(u)-1+H^1(u)),\\
H^1(u)&=x(G^1(u)+H^0(u)).
\end{array}\right.
\end{equation*}
Solving this system, we find:
$$F^0(u)=1+{\frac {ux \left( (u^2x^3+x^2) g_0^0+{u}^{2}{x}^{3}+{u}^{2}{x}^
{2}-{u}^{2}+{x}^{2} \right) }{{u}^{4}{x}^{2}+2\,{u}
^{2}{x}^{3}+{u}^{2}{x}^{2}-{u}^{2}+{x}^{2}}}
,
$$
$$F^1(u)=\frac{u^{2} x^{2} \left( x(\,x+1)g_0^0 -u^{2}+x^{2}\right)}{u^{4} x^{2}+2 u^{2} x^{3}+u^{2} x^{2}-u^{2}+x^{2}},\quad 
G^0(u)=-\frac{x^{2} \left( (\,u^{2} x^{2}-1)g_0^0+u^{2} x^{2}+u^{2} x  \right)}{u^{4} x^{2}+2 u^{2} x^{3}+u^{2} x^{2}-u^{2}+x^{2}},$$
$$G^1(u)=-\frac{\left(
(\,u^{2} x^{2}+x^{3}+x^2-1)g_0^0+u^{2} x^{2}+x^{3} \right) x u}{u^{4} x^{2}+2 u^{2} x^{3}+u^{2} x^{2}-u^{2}+x^{2}},$$
$$H^0(u)=\frac{u \,x^{2} \left(x(x+1)g_0^0 -u^{2}+x^{2}\right)}{u^{4} x^{2}+2 u^{2} x^{3}+u^{2} x^{2}-u^{2}+x^{2}}, \quad 
H^1(u)=-\frac{u \,x^{2} \left( (\,u^{2} x^{2}-1)g_0^0+u^{2} x^{2}+u^{2} x \right)}{u^{4} x^{2}+2 u^{2} x^{3}+u^{2} x^{2}-u^{2}+x^{2}}.
$$

In order to compute $g_0^0=G^0(0)$, we once again apply the kernel method. This consists in canceling the denominator of $F^1(u)$ by finding $u^2$ as an algebraic function $r$ of~$x$. If we substitute $u^2$ by $r$ in the numerator, the latter must vanish (to compensate for the  cancellation of the denominator), which determines  the value of $g_0^0$. We factorize the denominator $u^{4} x^{2}+2 u^{2} x^{3}+u^{2} x^{2}-u^{2}+x^{2}=x^2(u^2-r)(u^2-s)$, where
$$r=\frac{1-2 x^{3}-x^{2}-\sqrt{\Delta}}{2 x^{2}}, \quad s=\frac{1-2 x^{3}-x^{2}+\sqrt{\Delta}}{2 x^{2}},$$
and 
$\Delta=4 x^{6}+4 x^{5}-3 x^{4}-4 x^{3}-2 x^{2}+1.$
The cancellation of the numerator $(x^2 + x)g_0^0 - u^2 + x^2$ gives 
$$g^0_0=\frac{r-x^{2}}{x \left(x +1\right)}=\frac{1-2 x^{4}-2 x^{3}-x^{2}-\sqrt{\Delta}}{2 x^{3} \left(x +1\right)}.$$

Substituting this expression for $g_0^0$ into the generating functions above, and simplifying both the numerator and denominator by the factor $(u^2 - r)$, we obtain the following expressions:

\begin{thm}\label{teo1b} We have 
\[F^0(u)=1-\frac{u \left(x\,s +1\right)}{\left(x +1\right) \left(u^{2}-s \right)},
\quad  F^1(u)=\frac{-u^2}{u^2-s},
\quad G^0(u)=\frac{-2x^2-(1+r)x}{(1+x)(u^2-s)},  \]
\[G^1(u)=\frac{-u(r+x)}{(1+x)(u^2-s)},\quad H^0(u)=\frac{-u}{u^2-s},
\quad  \text{and} \quad  H^1(u)=\frac{u(-2x^2-(1+r)x)}{(1+x)(u^2-s)}.\]   
Finally, the bivariate generating function $S(x,u)$, where the coefficient of~$x^nu^k$ is the number of paths with $n$ steps ending at height $k$, satisfies
\begin{align*}S(x,u)&=\frac{u(-2x^2-x(r+2)-1)-sx-1}{x(1+x)(u^2-s)}.
\end{align*}
\end{thm}

The first terms of the series expansion are 
\begin{align*}1+ux+(u^2+u)x^2&+(u^3+2u+1)x^3+(u^4+u^3+u^2+3u+1)x^4+\\&(u^5+3u^3+3u^2+5u+1)x^5+
(u^6+u^5+2u^4+6u^3+3u^2+9u+3)x^6+\\&(u^7+4u^5+5u^4+11u^3+6u^2+15u+\bm{5})x^7 +O(x^8).
\end{align*}

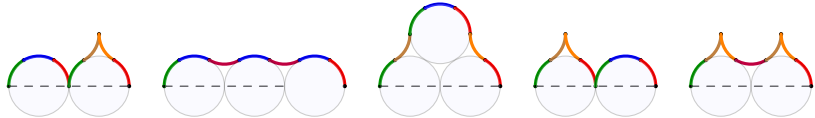
\begin{figure}[ht!]
  \centering
\begin{tikzpicture}[scale=0.4, line cap=round]
  \draw[dashed] (0,0)--(4,0);
  \renewcommand{\accX}{0}
  \renewcommand{\accY}{0}
  \po\NE\po\E\po\SE\po\NE\po\NEB\po\SEB\po\SE\po
  \filldraw[black] (0,0) circle (0.02);
  \filldraw[black] (\accX,\accY) circle (0.02);
  \draw[fill=blue!10,opacity=0.2] (1,0) circle (1);
  \draw[fill=blue!10,opacity=0.2] (3,0) circle (1);
\end{tikzpicture}\quad
\begin{tikzpicture}[scale=0.4, line cap=round]
  \draw[dashed] (0,0)--(4,0);
  \renewcommand{\accX}{0}
  \renewcommand{\accY}{0}
  \po\NE\po\E\po\EB\po\E\po\EB\po\E\po\SE\po
  \filldraw[black] (0,0) circle (0.02);
  \filldraw[black] (\accX,\accY) circle (0.02);
  \draw[fill=blue!10,opacity=0.2] (1,0) circle (1);
  \draw[fill=blue!10,opacity=0.2] (3,0) circle (1);
  \draw[fill=blue!10,opacity=0.2 ] (5,0) circle(1);
\end{tikzpicture}\quad
\begin{tikzpicture}[scale=0.4, line cap=round]
  \draw[dashed] (0,0)--(4,0);
  \renewcommand{\accX}{0}
  \renewcommand{\accY}{0}
  \po\NE\po\NEB\po\NE\po\E\po\SE\po\SEB\po\SE\po
  \filldraw[black] (0,0) circle (0.02);
  \filldraw[black] (\accX,\accY) circle (0.02);
  \draw[fill=blue!10,opacity=0.2] (1,0) circle (1);
  \draw[fill=blue!10,opacity=0.2] (3,0) circle (1);
  \draw[fill=blue!10,opacity=0.2 ] (2,1.75) circle(1);
\end{tikzpicture}\quad
\begin{tikzpicture}[scale=0.4, line cap=round]
  \draw[dashed] (0,0)--(4,0);
  \renewcommand{\accX}{0}
  \renewcommand{\accY}{0}
  \po\NE\po\NEB\po\SEB\po\SE\po\NE\po\E\po\SE\po
  \filldraw[black] (0,0) circle (0.02);
  \filldraw[black] (\accX,\accY) circle (0.02);
  \draw[fill=blue!10,opacity=0.2] (1,0) circle (1);
  \draw[fill=blue!10,opacity=0.2] (3,0) circle (1);
\end{tikzpicture}\quad
\begin{tikzpicture}[scale=0.4, line cap=round]
  \draw[dashed] (0,0)--(4,0);
  \renewcommand{\accX}{0}
  \renewcommand{\accY}{0}
  \po\NE\po\NEB\po\SEB\po\EB\po\NEB\po\SEB\po\SE\po
  \filldraw[black] (0,0) circle (0.02);
  \filldraw[black] (\accX,\accY) circle (0.02);
  \draw[fill=blue!10,opacity=0.2] (1,0) circle (1);
  \draw[fill=blue!10,opacity=0.2] (3,0) circle (1);
\end{tikzpicture}

\caption{The 5 paths with 7 steps (ending on the $x$-axis at abscissa $4$).}
\label{fig:pathsection2}
\end{figure}

In Figure~\ref{fig:pathsection2} we show the paths with 7 steps ending on the $x$-axis, corresponding to the bold coefficient in the above series.

Using the Vieta relation $rs=1$, we deduce the following result.

\begin{cor}\label{cor12} We have 
$$[u^{2k-1}]F^0(u)=\frac{sx+1}{x+1}r^{k}, ~k\geq 1,\quad  [u^{2k}]F^1(u)=r^{k}, ~k\geq 1,$$
$$[u^{2k}]G^0(u)=\frac{3x^2+r}{x(2x+r+1)}r^{k+1}, ~k\geq 1, \quad 
[u^{2k-1}]G^1(u)=\frac{x+r}{x+1}r^{k}, ~k\geq 1, \
$$ $$[u^{2k-1}]H^0(u)=r^k,~k\geq 1,\quad [u^{2k-1}]H^1(u)=\frac{3x^2+r}{x(2x+r+1)}r^{k},~k\geq 1,$$ $$[u^{2k}]S(x,u)=\frac{sx+1}{x(x+1)}r^{k+1}
,~k\geq 0, \qquad [u^{2k-1}]S(x,u)=\frac{rx+2x^2+2x+1}{x(x+1)}r^{k}
,~k\geq 1, \mbox{ and}$$ 
and all other coefficients are equal to zero, except for $[u^{0}]F^0(u) = 1$.

\end{cor}

In the following theorem we give an expression to compute the  coefficient $x^nu^k$ in the series expansion of $S(x,u)$. 

\begin{thm} Let $s(n,k)$ denote the number of paths with $n$ steps ending at height  $k$.  Then
\begin{align*} 
s(n,2k)&=w(n,k)+ w(n+1,k+1),\\
s(n,2k-1)&=w(n,k+1)+2w(n-1,k)+2w(n,k)+w(n+1,k),
\end{align*}
and $s(n,k)=0$ otherwise, where 
$$w(n,k)=\sum\limits_{i=0}^n (-1)^{i} \sum_{s=0}^{\left\lfloor\frac{n-i-2k}{4}\right\rfloor}\frac{k}{k+2s}\binom{k+2s}{s}t(n-i-4s-2k,2s+k),$$ and 
$$t(n,k)= \sum_{\ell=0}^{\lfloor n/2 \rfloor} 
\binom{k+\ell-1}{\ell}\binom{\ell}{n-2\ell}2^{n-2\ell}.$$
    \end{thm}
    
\begin{proof}
To obtain closed forms for the coefficients of $x^n u^k$, we first derive a closed form for the coefficient of $x^n$ in $r^k$, $k \geq 0$. Observe that  
\[
r(x) = \frac{1 - 2x^{3} - x^{2} - \sqrt{\Delta}}{2x^{2}}
= t(x)\, C(t^2(x)),
\]
where  
\[
C(x) = \frac{1 - \sqrt{1-4x}}{2x} \quad \text{and} \quad 
t(x) = \frac{x^2}{1 - x^2 - 2x^3}.
\]
Note that $C(x)$ is the generating function of the Catalan numbers. From \cite[equation 5.70]{GKP},  we have
\begin{align}
C^n(x)=\sum_{k=0}^\infty\frac{n}{n+2k} \binom{n+2k}{k}x^k. \label{powerCat}
\end{align}
 This implies
 \begin{align*}
   [x^n]r^m &= [x^n]\, t^m(x)\, C^m\!\bigl(t^2(x)\bigr) \\
   &= [x^n]\sum_{k \geq 0} \frac{m}{m+2k} \binom{m+2k}{k}\, t^{2k+m}(x) \\
   &= [x^n] \sum_{k \geq 0} \frac{m}{m+2k} \binom{m+2k}{k}\, 
      \frac{x^{4k+2m}}{(1 - x^2 - 2x^3)^{2k+m}}.
\end{align*}
Let   $t(n,m) = [x^n](1 - x^2 - 2x^3)^{-m}$. It is not difficult to verify that  
\begin{align}\label{eqt}
t(n,m) = \sum_{\ell=0}^{\lfloor n/2 \rfloor} 
\binom{m+\ell-1}{\ell}\binom{\ell}{n-2\ell}2^{n-2\ell}.\end{align}
Therefore, 
\begin{align*}
   [x^n]r^m &= [x^n]\sum_{k \geq 0}\sum_{\ell \geq 0} 
   \frac{m}{m+2k}\binom{m+2k}{k}\, t(\ell,2k+m)\, x^{4k+2m+\ell}.
\end{align*}
Setting $s = 4k+2m+\ell$, we obtain  
\begin{align*}
   [x^n]r^m&= [x^n]\sum_{k\geq 0}\sum_{s=4k+m} \frac{m}{m+2k}\binom{m+2k}{k}t(s-4k-2m,2k+m)x^{s}\\
   &=\sum_{k=0}^{\left\lfloor\frac{n-2m}{4}\right\rfloor}\frac{m}{m+2k}\binom{m+2k}{k}t(n-4k-2m,2k+m).
\end{align*}

Notice that we can also obtain the following identity directly from Lagrange inversion. 


 \[[x^n]r^m  = \sum_{\substack{\ell=0 \\ \ell-m \equiv 0 \pmod{2}}}^{\lfloor n/2 \rfloor}
      \frac{\ell}{n}\binom{\ell}{\tfrac{\ell-m}{2}} \, t(n-2\ell,\ell).\]
Using Corollary~\ref{cor12} and setting $u(n,k):=[x^{n}]r^k$, we can then obtain closed expressions for the numbers $f^0(n,k)$, $f^1(n,k)$,  $g^0(n,k)$, $g^1(n,k)$,  $h^0(n,k)$, and $h^1(n,k)$, which respectively count the partial packing paths of length $n$ ending at height $k$ with final steps $NE$,  $\overline{NE}$, $E$, $\overline{E}$, $SE$, and $\overline{SE}$.
\end{proof}

\begin{cor}
The generating function $S(x)=S(x,1)$, where  the coefficient of~$x^n$ is the  total number of paths with $n$ steps, satisfies
$$S(x)=\frac{1 - x - 2 x^2 -\sqrt{1 - 2 x + x^2 - 4 x^3 + 4 x^4}}{2 x^2 (-1 + 2 x)}.$$
An asymptotic approximation of the $n$-th term is
$$ \frac{2^{n+1}}{\sqrt{\pi n}}.$$
\end{cor}
The first terms of the series expansion are 
$$1+x+2x^2+4x^3+7x^4+13x^5+25x^6+47x^7+89x^8+171x^9+328x^{10}+630x^{11}+O(x^{12}).$$
This sequence does not appear in the OEIS.

\begin{cor} The generating function whose coefficient of~$x^n$ is the number of paths  with $n$ steps  ending on the $x$-axis satisfies
$$S(x,0)=\frac{1-x-\sqrt{(1 - x) (1 - 2 x) (1 + x + 2 x^2)}}{2 x^3}.$$
Moreover, we have $[x^0]S(x,0)=1$, $$[x^{n}]S(x,0)=\sum\limits_{k=0}^{\lfloor n/3\rfloor}\frac{1}{k+1}{2k\choose k}{n-2k-1\choose n-3k}, \mbox{ for } n\geq 0,$$   and 0 otherwise.  
An asymptotic approximation of the $n$-th term is 
$$ \frac{2^{n+1}}{\sqrt{\pi}n^{3/2}}.$$
\end{cor}

The first terms of the series expansion are 
$$1+x^3+x^4+x^5+3x^6+5x^7+7x^8+14x^9+26x^{10}+43x^{11}+79x^{12}+148x^{13}+O(x^{14}).$$
The sequence of coefficients of $x^{n}$ for $n \geq 0$ corresponds to \oeis{A346503}. At present, this sequence has no known combinatorial interpretation in the OEIS.

By calculating $\partial_u(S(x,u))\rvert_{u=1}$, and using classical methods~\cite{flajolet,orlov} for an asymptotic approximation of the coefficient of~$x^n$, we obtain the following.

\begin{cor} An asymptotic for the expected height in all paths of a given number of steps is given by $$\frac{17\sqrt{\pi n}}{27}.$$
\end{cor}

\subsection{A bijection with certain Motzkin paths.} 

Notice that $S(x,0)$ satisfies the functional equation $$S(x,0)=1+\frac{x^3}{1-x}S(x,0)^2.$$
This is also the functional equation of quarter-plane paths with $n$ steps, starting at the origin and ending on the $x$-axis, consisting of the steps $u=(1,1)$, $d=(1,-1)$ and $f=(1,0)$ , where 
down-step $d$ is always followed by at least one step $f$, and  each step $f$ follows either a step $d$ or another step $f$. Let $\mathcal{Q}_n$ denote the set of such paths with $n$ steps, and set $\mathcal{Q}=\cup_{n\geq 0}\mathcal{Q}_n$.

Proposition \ref{bij2prop} presents a constructive bijection $g$ between packing paths of length $n$ ending on the $x$-axis and paths in $\mathcal{Q}_n$. 

\begin{prop}\label{bij2prop} Define a  map $g:\mathcal{P}'\rightarrow\mathcal{Q}$ recursively as follows. For the empty path we set $g(\epsilon)=\epsilon$, and for $P\neq \epsilon$, we define 
$$g(P)=\left\{ \begin{array}{lll}
ug(Q) d ff \left(ff\right)^k g(R)&\mbox{ if } P= U \left(F\overline{F}\right)^k\overline{U}Q\overline{D}D R\quad &\mbox{ with } k\geq 0 \mbox{ and } Q,R\in \mathcal{P}',\\
ug(Q) d f \left(ff\right)^k g(R)&\mbox{ if } P= U \left(F\overline{F}\right)^k\overline{U}Q\overline{D}F S\quad &\mbox{ with } k\geq 0 \mbox{ and } Q,US\in \mathcal{P}',\\
ug(R) d f \left(ff\right)^k &\mbox{ if } P= U F\left(\overline{F}F\right)^kD R\quad &\mbox{ with } k\geq 0 \mbox{ and } Q,R\in \mathcal{P}',
\end{array}
\right.$$
Then $g$ is a bijection from $\mathcal{P}'$ onto $\mathcal{Q}$.
\end{prop}

Due to the recursive definition, the image by $g$ of a packing path in $\mathcal{P}'_{n}$ is path in $\mathcal{Q}_n$, and it is easy to see that $g$ is a bijection. We refer to Figure~\ref{biection2} for an illustration of $g$ when $P=U\overline{U}~UFD~\overline{D}D\overline{F}U\overline{D}DUF\overline{F}FD$ and
$f(P)=uudfdfudffudfff$.

\begin{figure}[H]
    \centering
\begin{tikzpicture}[scale=0.4, line cap=round]
\draw[dashed] (0,0)--(10,0);
\renewcommand{\accX}{0}
\renewcommand{\accY}{0}
\po\NE\po\NEB\po\NE\po
\E\po\SE\po\SEB\po
\EB\po\NEB\po\SEB\po\SE
\po\NE\po
\E\po\EB\po\E\SE\po
\filldraw[black] (0,0) circle (0.02);
\filldraw[black] (\accX,\accY) circle (0.02);
\draw[fill=blue!10,opacity=0.2 ] (1,0) circle(1);
\draw[fill=blue!10,opacity=0.2 ] (2,1.75) circle(1);
\draw[fill=blue!10,opacity=0.2 ] (3,0) circle(1);
\draw[fill=blue!10,opacity=0.2 ] (5,0) circle(1);
\draw[fill=blue!10,opacity=0.2 ] (7,0) circle(1);
\draw[fill=blue!10,opacity=0.2 ] (9,0) circle(1);
\end{tikzpicture}\qquad $\longrightarrow$\qquad
\begin{tikzpicture}[scale=0.4, line cap=round]
\draw[dashed] (0,0)--(15,0);
\draw[line width=0.7pt] (0,0)--(1,1)--(2,2)--(3,1)--(4,1)--(5,0)--(6,0)--(7,1)--(8,0)--(9,0)--(10,0)--(11,1)--(12,0)--(13,0)--(14,0)--(15,0);
\point{0}{0}
\point{1cm}{1cm}
\point{2cm}{2cm}\point{3cm}{1cm}\point{4cm}{1cm}\point{5cm}{0cm}\point{6cm}{0cm}\point{7cm}{1cm}\point{8cm}{0cm}
\point{9cm}{0cm}\point{10cm}{0cm}\point{11cm}{1cm}\point{12cm}{0cm}\point{13cm}{0cm}\point{14cm}{0cm}\point{15cm}{0cm}
\end{tikzpicture}
  \caption{The packing path $P=U\overline{U}~UFD~\overline{D}D\overline{F}U\overline{D}DUF\overline{F}FD$ and its image by $g$: 
$f(P)=uudfdfudffudfff$.}
    \label{biection2}
\end{figure}

\subsection{A connection with Riordan arrays}
Consider the matrices $\mathcal{S}_e=[s(n+k,2k)]_{n, k\geq 0}$ and   $\mathcal{S}_o=[s(n+k+1,2k+1)]_{n, k\geq 0}$.  The first few rows of $\mathcal{S}_e$ and $\mathcal{S}_o$ are
\[
\footnotesize 
\mathcal{S}_e=\left(
\begin{array}{ccccccc}
 1 & 0 & 0 & 0 & 0 & 0 & 0 \\
 0 & 1 & 0 & 0 & 0 & 0 & 0 \\
 0 & 0 & 1 & 0 & 0 & 0 & 0 \\
 1 & 1 & 0 & 1 & 0 & 0 & 0 \\
 1 & 3 & 2 & 0 & 1 & 0 & 0 \\
 1 & 3 & 5 & 3 & 0 & 1 & 0 \\
 3 & 6 & 6 & 7 & 4 & 0 & 1 \\
\end{array} 
\right), 
\quad 
\mathcal{S}_o=\left(
\begin{array}{cccccccc}
 1 & 0 & 0 & 0 & 0 & 0 & 0 \\
 1 & 1 & 0 & 0 & 0 & 0 & 0 \\
 2 & 1 & 1 & 0 & 0 & 0 & 0 \\
 3 & 3 & 1 & 1 & 0 & 0 & 0 \\
 5 & 6 & 4 & 1 & 1 & 0 & 0 \\
 9 & 11 & 9 & 5 & 1 & 1 & 0 \\
 15 & 22 & 18 & 12 & 6 & 1 & 1 \\
\end{array}
\right).
\]

Both matrices are Riordan arrays.  From Corollary~\ref{cor12} and the definition of a Riordan array, we obtain the following result.  

\begin{thm}\label{Riordan2}
The matrices $\Sr_e$ and $\Sr_o$ are Riordan arrays. More precisely,
$$\Sr_e=\left(\frac{1 - x -\sqrt{p(x)}}{2x^3}, \frac{1 - x^2 (1 + 2 x) - (1 + x)\sqrt{p(x)}}{2x^3} \right)$$
and
$$\Sr_o=\left(\frac{1 - 3 x^3 - 2 x^4 - (1 + x + x^2)\sqrt{p(x)}}{2x^5}, \frac{1 - x^2 - 2 x^3 - (1 + x)\sqrt{p(x)}}{2x^3} \right),$$
where $p(x)=1 - 2 x + x^2 - 4 x^3 + 4 x^4$.
\end{thm}

\section{Packing paths with a given area}

In this section, we enumerate packing paths ending on the $x$-axis according to their area and width.  Recall that the area of a packing path is defined as the number of circles lying below the path and above the line $y=-2$, while the width is the abscissa of its endpoint.

Let $\mathcal{A}_n$ denote the set of packing paths  that start at the origin, end on the $x$-axis at abscissa $n$, and remain in the quarter plane.  Let $\mathcal{A} = \bigcup_{n \geq 0} \mathcal{A}_n$. Define the bivariate generating function
\[A(x,y)=\sum_{P\in \mathcal{A}}x^{\s(P)}y^{\ar(P)},\]
where $x$ marks the width and $y$ marks the area.

We also consider the set $\mathcal{B}$ of paths in the quarter plane that start at the point $(1, \sqrt{3})$, end with a step of type $F$ or $\overline{D}$ on the line $y = \sqrt{3}$, and do not touch the $x$-axis. The width of a path in $\mathcal{B}$ is defined as the abscissa of its endpoint minus one (to compensate for the fact that the path starts at abscissa one). Let $\mathcal{B}_n$ denote the set of paths in $\mathcal{B}$ of width $n$. The area of a path in $\mathcal{B}$ is defined as the number of full  circles below the path and strictly above the $x$-axis. 

For instance, the area (resp. the width) of the path $P=F\overline{F}F\overline{F}F$ is 0 (resp. 10), and the area (resp. width) of the path $Q=\overline{U}UFD\overline{D}$ is 1 (resp. 6).

Let $B(x, y)$ be the bivariate generating function that counts the number of paths in $\mathcal{B}$, where $x$ marks the width and $y$ marks the area. 

By means of a recursive decomposition on the first return to the $x$-axis, we obtain the following continued fraction representation for the generating function.

\begin{thm} The generating function for the number of packing paths with respect to their  width and area is given by \[
A(x,y)
=\frac{1}{1-
\cfrac{2x^{4}y}{1-2x^{4}y-\cfrac{x^{4}y^{2}}{1-3x^{4}y^{2}-\cfrac{x^{4}y^{3}}{1-3x^{4}y^{3}-\ddots}}}
\;}.
\]
\end{thm}

\begin{proof}
Any path $P\in\mathcal{A}$ can be decomposed with respect to its first return  on the $x$-axis as  $P=U Q D R $, where $Q\in \mathcal{B}$ and $R\in\mathcal{A}$. 
The area of $P$ is the sum of the areas of   $UQD$  and $R$. 

The area of \( UQD \), in turn, is equal to the area of \( Q\in\mathcal{B} \) plus the number of  full circles lying between \( Q \) and the line $y=-2$,   that is,
\[
\frac{2 + \texttt{width}(Q)}{4}.
\]
Hence 
$$\ar(UQDR)=\frac{2 + \texttt{width}(Q)}{4}+\ar(R).$$
Taking this recursive decomposition into account, we derive the following functional equation for \( A(x, y) \):
\[
    A(x,y)=1+x^2y^{1/2}B(xy^{1/4},y)A(x,y),\]
and therefore
\begin{equation}\label{ecCF1}
     A(x,y)=\frac{1}{1-x^2y^{1/2}B(xy^{1/4},y)}.
\end{equation}

Now consider any path $P'\in\mathcal{B}$. According to its first return to the line  $y=\sqrt{3}$, it has one of  the following four forms: 
\begin{itemize}
   \item[($i$)] $P'=F$,  
   \item[($ii$)] $P'=F\overline{F} R$  with  $R\in \mathcal{B}$, 
    \item[($iii$)] $P'=\overline{U} Q \overline{D}$  with $Q\in\mathcal{A}$,  
    \item[($iv$)] $P'=\overline{U} Q \overline{D}\overline{F}R$  with $Q \in\mathcal{A}$ and  $R\in\mathcal{B}$.
\end{itemize}
 The contribution for the path satisfying ($i$) is $x^2$; The contribution for the paths satisfying ($ii$) is $x^4B(x,y)$; The contribution for the paths satisfying ($iii$) is $x^2A(x,y)$;  The contribution for the paths satisfying ($iv$) is $x^4A(x,y)B(x,y)$.
Summing these cases yields the functional equation
\begin{equation}\label{ecCF2}
B(x,y)=x^2+x^4B(x,y)+x^2A(x,y)+x^4A(x,y)B(x,y),
\end{equation}
Combining \eqref{ecCF1} and \eqref{ecCF2}, we obtain
\begin{equation}
B(x,y)=x^2+x^4B(x,y)+\frac{x^2(1+x^2B(x,y))}{1-x^2y^{1/2}B(xy^{1/4},y)}.
    \label{funcarea}
\end{equation}
From \eqref{funcarea}, it follows that
\[
B(x,y)=\frac{x^2\bigl(2 - x^2 y^{1/2}B(x y^{1/4},y)\bigr)}
{1-2x^4-(1-x^4)x^2 y^{1/2}B(x y^{1/4},y)}.
\]
Let \(B_k:=B(xy^{k/4},y)\) and \(t_k:=xy^{k/4}\). Define the auxiliary quantities
\[
 U_k:=1-t_k^{4},\qquad V_k:=1-2t_k^{4}, \qquad W_k:=t_k^{2}y^{1/2}.
\]
Then, in particular,
\[
B(x,y)=B_0=\frac{t_0^{2}\bigl(2-W_0B_1\bigr)}{V_0-U_0W_0B_1}.
\]
In general,
\[
B_k=\frac{t_k^{2}\bigl(2-W_kB_{k+1}\bigr)}{V_k-U_k W_kB_{k+1}},\qquad k\geq 0.
\]
Since $V_k+1=2U_k$,  we can rewrite this as
\begin{multline*}
   B_k= \frac{t_k^2(2-W_kB_{k+1})}{V_k-U_kW_kB_{k+1}}=\frac{2t_k^2(2-W_kB_{k+1})}{2V_k-2U_kW_kB_{k+1}}=\frac{2t_k^2(2-W_kB_{k+1})}{2V_k-(V_k+1)W_kB_{k+1}}\\=\frac{2t_k^2(2-W_kB_{k+1})}{V_k(2-W_kB_{k+1})-W_kB_{k+1}} = \frac{2t_k^2}{V_k-\frac{W_kB_{k+1}}{2-W_kB_{k+1}}}=\frac{2t_k^2}{V_k-H_k},
\end{multline*}
where 
$$H_k:=\frac{W_kB_{k+1}}{2-W_kB_{k+1}}.$$
Substituting \(B_{k+1}\) into \(H_k\) gives
$$H_{k}=\frac{2t_{k+1}^2W_k}{2(V_{k+1}-H_{k+1})-2t_{k+1}^2W_{k}}=\frac{t_{k+1}^2W_k}{V_{k+1}-t_{k+1}^2W_{k}-H_{k+1}}.$$
Therefore,
$$B_k=\cfrac{2t_k^2}{V_k-\cfrac{t_{k+1}^2W_k}{V_{k+1}-t_{k+1}^2W_{k}-H_{k+1}}}.$$
Iterating this identity and setting \(k=0\) yields
\begin{align*}
 B_0=\cfrac{2t_0^2}{V_0-\cfrac{t_{1}^2W_0}{V_{1}-t_{1}^2W_{0}-\cfrac{t_2^2W_1}{V_2-t_2^2W_1-\cfrac{t_{3}^2W_2}{V_{3}-t_{3}^2W_{2}-\ddots}}}}.
\end{align*}
Finally, replacing the auxiliary variables, we obtain the continued fraction:
\[
B(x,y)=
\cfrac{2x^{2}}{1-2x^{4}-\cfrac{x^{4}y}{
1-3x^{4}y-\cfrac{x^{4}y^{2}}{
1-3x^{4}y^{2}-\cfrac{x^{4}y^{3}}{
1-3x^{4}y^{3}-\ddots}}}}
.\]
Substituting this expression into \eqref{ecCF1} yields the desired form of \(A(x,y)\).
\end{proof}

The first few terms of the power series expansion are
\begin{align*}A(x,y)&=1+2 x^4 y+x^8 \left(2 y^3+8 y^2\right)+x^{12} \left(2 y^6 + \bm{8 y^5}+16 y^4+32
   y^3\right)\\&+x^{16} \left(2 y^{10}+8 y^9+16 y^8+48 y^7+72 y^6+96 y^5+128
   y^4\right)+O\left(x^{20}\right).
   \end{align*}

\begin{figure}[H]
    \centering
\begin{tikzpicture}[scale=0.4, line cap=round]
\draw[dashed] (0,0)--(6,0);
\renewcommand{\accX}{0}
\renewcommand{\accY}{0}
\po\NE\po\NEB\po\NE\po
\E\po\SE\po\NE\po\E\po
\SE\po\SEB\po\SE\po
\filldraw[black] (0,0) circle (0.02);
\filldraw[black] (\accX,\accY) circle (0.02);\draw[fill=blue!10,opacity=0.2 ] (1,0) circle(1);\draw[fill=blue!10,opacity=0.2 ] (3,0) circle(1);\draw[fill=blue!10,opacity=0.2 ] (5,0) circle(1);
\draw[fill=blue!10,opacity=0.2 ] (2,1.75) circle(1);\draw[fill=blue!10,opacity=0.2 ] (4,1.75) circle(1);
\end{tikzpicture}\quad
\begin{tikzpicture}[scale=0.4, line cap=round]
\draw[dashed] (0,0)--(6,0);
\renewcommand{\accX}{0}
\renewcommand{\accY}{0}
\po\NE\po\NEB\po\NE\po
\E\po\SE\po\NE\po\NEB\po\SEB\po
\SE\po\SEB\po\SE\po
\filldraw[black] (0,0) circle (0.02);
\filldraw[black] (\accX,\accY) circle (0.02);\draw[fill=blue!10,opacity=0.2 ] (1,0) circle(1);\draw[fill=blue!10,opacity=0.2 ] (3,0) circle(1);\draw[fill=blue!10,opacity=0.2 ] (5,0) circle(1);
\draw[fill=blue!10,opacity=0.2 ] (2,1.75) circle(1);\draw[fill=blue!10,opacity=0.2 ] (4,1.75) circle(1);
\end{tikzpicture}\quad
\begin{tikzpicture}[scale=0.4, line cap=round]
\draw[dashed] (0,0)--(6,0);
\renewcommand{\accX}{0}
\renewcommand{\accY}{0}
\po\NE\po\NEB\po\NE\po
\E\po\EB\po\E\po
\SE\po\SEB\po\SE\po
\filldraw[black] (0,0) circle (0.02);
\filldraw[black] (\accX,\accY) circle (0.02);\draw[fill=blue!10,opacity=0.2 ] (1,0) circle(1);\draw[fill=blue!10,opacity=0.2 ] (3,0) circle(1);\draw[fill=blue!10,opacity=0.2 ] (5,0) circle(1);
\draw[fill=blue!10,opacity=0.2 ] (2,1.75) circle(1);\draw[fill=blue!10,opacity=0.2 ] (4,1.75) circle(1);
\end{tikzpicture}\quad
\begin{tikzpicture}[scale=0.4, line cap=round]
\draw[dashed] (0,0)--(6,0);
\renewcommand{\accX}{0}
\renewcommand{\accY}{0}
\po\NE\po\NEB\po\NE\po
\E\po\EB\po\NEB\po\SEB\po
\SE\po\SEB\po\SE\po
\filldraw[black] (0,0) circle (0.02);
\filldraw[black] (\accX,\accY) circle (0.02);\draw[fill=blue!10,opacity=0.2 ] (1,0) circle(1);\draw[fill=blue!10,opacity=0.2 ] (3,0) circle(1);\draw[fill=blue!10,opacity=0.2 ] (5,0) circle(1);
\draw[fill=blue!10,opacity=0.2 ] (2,1.75) circle(1);\draw[fill=blue!10,opacity=0.2 ] (4,1.75) circle(1);
\end{tikzpicture}\\[0.1cm]
\begin{tikzpicture}[scale=0.4, line cap=round]
\draw[dashed] (0,0)--(6,0);
\renewcommand{\accX}{0}
\renewcommand{\accY}{0}
\po\NE\po\NEB\po\NE\po
\po\NEB\po\SEB\po\SE\po\NE\po\E\po
\SE\po\SEB\po\SE\po
\filldraw[black] (0,0) circle (0.02);
\filldraw[black] (\accX,\accY) circle (0.02);\draw[fill=blue!10,opacity=0.2 ] (1,0) circle(1);\draw[fill=blue!10,opacity=0.2 ] (3,0) circle(1);\draw[fill=blue!10,opacity=0.2 ] (5,0) circle(1);
\draw[fill=blue!10,opacity=0.2 ] (2,1.75) circle(1);\draw[fill=blue!10,opacity=0.2 ] (4,1.75) circle(1);
\end{tikzpicture}\quad
\begin{tikzpicture}[scale=0.4, line cap=round]
\draw[dashed] (0,0)--(6,0);
\renewcommand{\accX}{0}
\renewcommand{\accY}{0}
\po\NE\po\NEB\po\NE\po
\po\NEB\po\SEB\po\SE\po\NE\po\NEB\po\SEB\po
\SE\po\SEB\po\SE\po
\filldraw[black] (0,0) circle (0.02);
\filldraw[black] (\accX,\accY) circle (0.02);\draw[fill=blue!10,opacity=0.2 ] (1,0) circle(1);\draw[fill=blue!10,opacity=0.2 ] (3,0) circle(1);\draw[fill=blue!10,opacity=0.2 ] (5,0) circle(1);
\draw[fill=blue!10,opacity=0.2 ] (2,1.75) circle(1);\draw[fill=blue!10,opacity=0.2 ] (4,1.75) circle(1);
\end{tikzpicture}\quad
\begin{tikzpicture}[scale=0.4, line cap=round]
\draw[dashed] (0,0)--(6,0);
\renewcommand{\accX}{0}
\renewcommand{\accY}{0}
\po\NE\po\NEB\po\NE\po
\po\NEB\po\SEB\po\EB\po\E\po
\SE\po\SEB\po\SE\po
\filldraw[black] (0,0) circle (0.02);
\filldraw[black] (\accX,\accY) circle (0.02);
\draw[fill=blue!10,opacity=0.2 ] (1,0) circle(1);\draw[fill=blue!10,opacity=0.2 ] (3,0) circle(1);\draw[fill=blue!10,opacity=0.2 ] (5,0) circle(1);
\draw[fill=blue!10,opacity=0.2 ] (2,1.75) circle(1);\draw[fill=blue!10,opacity=0.2 ] (4,1.75) circle(1);
\end{tikzpicture}\quad
\begin{tikzpicture}[scale=0.4, line cap=round]
\draw[dashed] (0,0)--(6,0);
\renewcommand{\accX}{0}
\renewcommand{\accY}{0}
\po\NE\po\NEB\po\NE\po
\po\NEB\po\SEB\po\EB\po\NEB\po\SEB\po
\SE\po\SEB\po\SE\po
\filldraw[black] (0,0) circle (0.02);
\filldraw[black] (\accX,\accY) circle (0.02);
\draw[fill=blue!10,opacity=0.2 ] (1,0) circle(1);\draw[fill=blue!10,opacity=0.2 ] (3,0) circle(1);\draw[fill=blue!10,opacity=0.2 ] (5,0) circle(1);
\draw[fill=blue!10,opacity=0.2 ] (2,1.75) circle(1);\draw[fill=blue!10,opacity=0.2 ] (4,1.75) circle(1);
\end{tikzpicture}
    \caption{The 8 paths ending on the $x$-axis at abscissa 12 above five circles}
    \label{fig:patharea}
\end{figure}
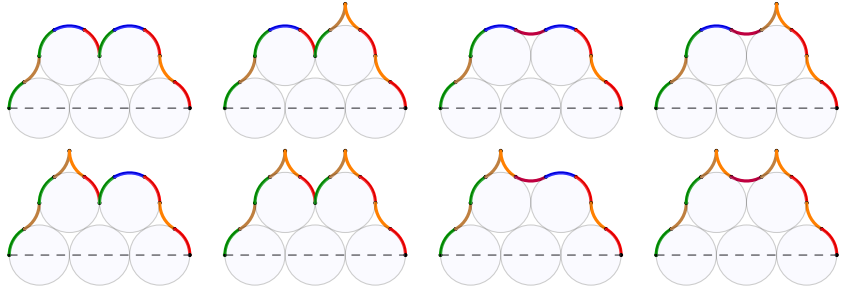

Figure~\ref{fig:patharea} shows the  packing paths of area $5$ and width $12$, corresponding to the bold coefficient in the above series.

Note that $A(x,1)$ is a reformulation of $g_0^0$, presented in Section 2, in term of a continued fraction.

\begin{cor} The generating function for the number of packing paths with respect to the area is given by 
$$A(1,y)=\frac{1}{1-
\cfrac{2y}{1-2y-\cfrac{y^{2}}{1-3y^{2}-\cfrac{y^{3}}{1-3y^{3}-\ddots}}}
\;}.$$
\end{cor}

The first terms of the series expansion of $A(1,y)$ are
$$1+2y+8y^2+34y^3+144y^4+616y^5+2634y^6+11280y^7+48304y^8+206904y^9+O(y^{10}).$$

\section{Packing paths with a given kissing number}

As previously defined, let $\mathcal{A}_n$ be the set of paths in the quarter plane that start at the origin, end on the $x$-axis at abscissa $n$, and let $\mathcal{A} = \bigcup_{n \geq 0} \mathcal{A}_n$. Similarly, let $\mathcal{B}_n$ denote the set of paths of width $n$ in the quarter plane that start at the point $(1, \sqrt{3})$, end with a step of type $F$ or $\overline{D}$ on the line $y = \sqrt{3}$, and do not touch the $x$-axis. We set $\mathcal{B} = \bigcup_{n \geq 1} \mathcal{B}_n$.

The kissing number of a packing path in $\mathcal{A}$ (resp. $\mathcal{B}$) is defined as the number of circles below the path that touch the path. In this section, we enumerate packing paths that end on the $x$-axis, according to  their width and the kissing number. 

Let $A(x, y)$ (resp. $B(x, y)$) denote the bivariate generating function that counts the number of paths in $\mathcal{A}$ (resp. $\mathcal{B}$), where $x$ marks the width and $y$ marks the kissing number.

\begin{thm} The generating function for the number of packing paths with respect to the width and the kissing number is given by \[
A(x,y)
={\frac {1 - 4 x^4 y + 3 x^4 y^2 -\sqrt {1 - 8 x^4 y - 2 x^4 y^2 + 16 x^8 y^2 - 8 x^8 y^3 + x^8 y^4}}{
4{x}^{4}{y}^{2}}}.
\]
\end{thm}

\begin{proof}
Any path $P$ in $\mathcal{A}$ can be decomposed with respect to its first return  on the $x$-axis as $P=U Q D R $, where $Q\in \mathcal{B}$ and $R\in\mathcal{A}$.  The kissing number of $P$ is the sum of the kissing number of $Q$ and $R$. Hence, we obtain the functional equation:
\begin{align}\label{ec1gfkiss}
    A(x,y) = 1+x^2B(x,y)A(x,y).
\end{align}
Similarly, any path $P'$ in $\mathcal{B}$ can be decomposed with respect to its first return  on the line $y=\sqrt{3}$ in the following four forms: 
\begin{itemize}
   \item[($i$)] $P'=F$,   
   \item[($ii$)] $P'=F\overline{F} R$  with  $R\in \mathcal{B}$,
    \item[($iii$)] $P'=\overline{U} Q \overline{D}$  with $Q\in\mathcal{A}$,
    \item[($iv$)] $P'=\overline{U} Q \overline{D}\overline{F}R$  with $Q \in\mathcal{A}$ and  $R\in\mathcal{B}$.
\end{itemize}
 The corresponding contributions for the path satisfying ($i$) is $x^2y$; The contribution for the paths satisfying  ($ii$) is $x^4yB(x,y)$; The contribution for the paths satisfying  ($iii$) is $x^2y$ if $Q=\epsilon$, and $x^2y^2(A(x,y)-1)$ otherwise; The contribution for the paths satisfying  ($iv$) is $x^4yB(x,y)$ if $Q=\epsilon$, and $x^4y^2B(x,y)(A(x,y)-1)$ otherwise. Considering the four cases, we obtain the
following functional equation
\begin{multline}\label{ec1gfkiss2}
    B(x,y) = x^4yB(x,y)+
    x^2y+x^2y+x^2y^2(A(x,y)-1)+\\
    x^4yB(x,y)+x^4y^2B(x,y)(A(x,y)-1).
\end{multline}
Combining \eqref{ec1gfkiss} and \eqref{ec1gfkiss2}  we obtain the desired result. 
\end{proof}

The first few terms of the power series expansion are
\begin{align*}A(x,y)=1+2x^4y&+(2y^3+8y^2)x^8+(2y^5+24y^4+32y^3)x^{12}+(2y^7+48y^6+192y^5+128y^4)x^{16}+\\
&(2y^9+80y^8+640y^7+1280y^6+512y^5)x^{20}+\\
&\quad\quad(2 y^{11}+120 y^{10}+1600 y^9+6400 y^8+7680 y^7+2048 y^6) x^{24} + O(x^{28}).
   \end{align*}

\begin{figure}[H]
    \centering
\begin{tikzpicture}[scale=0.4, line cap=round]
\draw[dashed] (0,0)--(12,0);
\renewcommand{\accX}{0}
\renewcommand{\accY}{0}
\po\NE\po\NEB\po\NE\po
\E\po\EB\po\NEB\po\NE\po
\NEB\po\NE\po\NEB\po\SEB
\po\EB\po
\E\po\SE\po\SEB\po\SE\po\SEB\po\SE\po\SEB\po\SE\po
\filldraw[black] (0,0) circle (0.02);
\filldraw[black] (\accX,\accY) circle (0.02);
\draw[fill=blue!10,opacity=0.2 ] (1,0) circle(1);
\draw[fill=blue!10,opacity=0.2 ] (2,1.75) circle(1);
\draw[fill=blue!10,opacity=0.2 ] (4,1.75) circle(1);
\draw[fill=blue!10,opacity=0.2 ] (5,3.5) circle(1);
\draw[fill=blue!10,opacity=0.2 ] (6,5.25) circle(1);
\draw[fill=blue!10,opacity=0.2 ] (8,5.25) circle(1);
\draw[fill=blue!10,opacity=0.2 ] (9,3.5) circle(1);
\draw[fill=blue!10,opacity=0.2 ] (10,1.75) circle(1);
\draw[fill=blue!10,opacity=0.2 ] (11,0) circle(1);
\end{tikzpicture}\quad
\begin{tikzpicture}[scale=0.4, line cap=round]
\draw[dashed] (0,0)--(12,0);
\renewcommand{\accX}{0}
\renewcommand{\accY}{0}
\po\NE\po\NEB\po\NE\po
\NEB\po\SEB\po\SE\po\SEB\po\EB\po\NEB\po\NE\po
\E\po\SE\po\SEB\po\SE\po\NE\po\NEB\po\NE\po
\E\po\SE\po\SEB\po\SE
\filldraw[black] (0,0) circle (0.02);
\filldraw[black] (\accX,\accY) circle (0.02);
\draw[fill=blue!10,opacity=0.2 ] (1,0) circle(1);
\draw[fill=blue!10,opacity=0.2 ] (2,1.75) circle(1);
\draw[fill=blue!10,opacity=0.2 ] (3,0) circle(1);
\draw[fill=blue!10,opacity=0.2 ] (5,0) circle(1);
\draw[fill=blue!10,opacity=0.2 ] (6,1.75) circle(1);
\draw[fill=blue!10,opacity=0.2 ] (7,0) circle(1);
\draw[fill=blue!10,opacity=0.2 ] (9,0) circle(1);
\draw[fill=blue!10,opacity=0.2 ] (10,1.75) circle(1);
\draw[fill=blue!10,opacity=0.2 ] (11,0) circle(1);
\end{tikzpicture}

  \caption{Two paths  ending on the $x$-axis at abscissa 24 and with a kissing number equal to 9 (there are 1600 such paths).}
    \label{fig:pathkiss}
\end{figure}
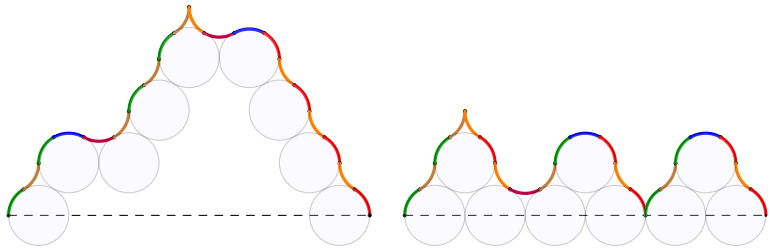

Figure \ref{fig:pathkiss} shows two paths with kissing number equal to 9 and width 24.

\begin{cor} The generating function for the number of packing paths with respect to the kissing number is given by 
$$A(1,y)={\frac {1 - 4 y + 3 y^2 - (1-y)\sqrt{1 - 6 y + y^2}}{4{y}^{2}}}.$$
\end{cor}

The first terms of the series expansion of $A(1,y)$ are
$$1+2y+8y^2+34y^3+152y^4+706y^5+3376y^6+16514y^7+82256y^8 + O(y^9).$$
The sequence of coefficients corresponds to twice the sequence \oeis{A239204}.  Let $r_n$ denote the number of large Schröder paths of length $n$, that is, quarter-plane paths with $n$ steps, starting and ending on the $x$ axis, consisting of the steps $u=(1,1), d=(1,-1)$, and $h=(2,0)$. It is well known that their generating function is
\[R(y)=\sum_{n\geq 0}r_ny^n=\frac{1-y-\sqrt{1-6y+y^2}}{2y}.\]
Therefore, we have the relation
$$A(1,y)=\frac{1-y}{y}(R(y)-1).$$
Comparing coefficients, we obtain
$$[y^n]A(1,y)=\frac{r_{n+1}-r_n}{2}, \qquad n\geq 1.$$
Using the known formula
\[
r_n=\sum_{k=0}^{n}\binom{n+k}{2k}C_k,
\]
where $C_k$ denotes the $k$-th Catalan number, we find that 
\[
[y^n]A(1,y) = \frac{1}{2}\sum_{k=1}^{n+1} \binom{n+k}{2k-1}C_k, \qquad n\geq 1.
\]
We leave as an open question the task of finding a direct combinatorial interpretation of this relation between the coefficients of $A(1,y)$ and the large Schröder numbers.

\subsection{Inchworm paths}
\begin{figure}[H]
    \centering
\begin{tikzpicture}[scale=0.4, line cap=round]
\draw[dashed] (0,0)--(12,0);
\renewcommand{\accX}{0}
\renewcommand{\accY}{0}
\po\NE
\po\E\po\EB\po
\po\E\po\EB\po
\po\E\po\EB\po
\po\E\po\EB\po
\po\E\po\EB\po
\po\E\po\SE\po
\draw[fill=blue!10,opacity=0.2 ] (1,0) circle(1);
\draw[fill=blue!10,opacity=0.2 ] (3,0) circle(1);
\draw[fill=blue!10,opacity=0.2 ] (5,0) circle(1);
\draw[fill=blue!10,opacity=0.2 ] (7,0) circle(1);
\draw[fill=blue!10,opacity=0.2 ] (9,0) circle(1);
\draw[fill=blue!10,opacity=0.2 ] (11,0) circle(1);
\end{tikzpicture}\quad
\begin{tikzpicture}[scale=0.4, line cap=round]
\draw[dashed] (0,0)--(10,0);
\renewcommand{\accX}{0}
\renewcommand{\accY}{0}
\po\NE
\po\E\po\EB\po
\po\E\po\EB\po
\po\E\po\EB\po
\po\NEB\po\NE
\po
\po\E\po\SE\po
\SEB\po\SE\po

\draw[fill=blue!10,opacity=0.2 ] (1,0) circle(1);
\draw[fill=blue!10,opacity=0.2 ] (3,0) circle(1);
\draw[fill=blue!10,opacity=0.2 ] (5,0) circle(1);
\draw[fill=blue!10,opacity=0.2 ] (7,0) circle(1);
\draw[fill=blue!10,opacity=0.2 ] (9,0) circle(1);
\draw[fill=blue!10,opacity=0.2 ] (8,1.75) circle(1);
\end{tikzpicture}
\quad
\begin{tikzpicture}[scale=0.4, line cap=round]
\draw[dashed] (0,0)--(10,0);
\renewcommand{\accX}{0}
\renewcommand{\accY}{0}
\po\NE
\po\E\po\EB\po
\po\E\po\EB\po
\po\NEB\po\NE
\po
\E\po\SE\po
\SEB\po\EB\po\E
\po\SE\po

\draw[fill=blue!10,opacity=0.2 ] (1,0) circle(1);
\draw[fill=blue!10,opacity=0.2 ] (3,0) circle(1);
\draw[fill=blue!10,opacity=0.2 ] (5,0) circle(1);
\draw[fill=blue!10,opacity=0.2 ] (7,0) circle(1);
\draw[fill=blue!10,opacity=0.2 ] (9,0) circle(1);
\draw[fill=blue!10,opacity=0.2 ] (6,1.75) circle(1);
\end{tikzpicture}
\quad
\begin{tikzpicture}[scale=0.4, line cap=round]
\draw[dashed] (0,0)--(10,0);
\renewcommand{\accX}{0}
\renewcommand{\accY}{0}
\po\NE
\po\E\po\EB\po

\NEB\po\NE
\po
\E\po\SE\po
\SEB
\po\EB\po\E
\po\EB\po\E
\po\SE\po

\draw[fill=blue!10,opacity=0.2 ] (1,0) circle(1);
\draw[fill=blue!10,opacity=0.2 ] (3,0) circle(1);
\draw[fill=blue!10,opacity=0.2 ] (5,0) circle(1);
\draw[fill=blue!10,opacity=0.2 ] (7,0) circle(1);
\draw[fill=blue!10,opacity=0.2 ] (9,0) circle(1);
\draw[fill=blue!10,opacity=0.2 ] (4,1.75) circle(1);
\end{tikzpicture}
\quad
\begin{tikzpicture}[scale=0.4, line cap=round]
\draw[dashed] (0,0)--(10,0);
\renewcommand{\accX}{0}
\renewcommand{\accY}{0}
\po\NE\po\NEB\po\NE\po\E\po\SE\po\SEB\po\EB\po

\po\E\po\EB\po\E
\po\EB\po\E
\po\SE\po

\draw[fill=blue!10,opacity=0.2 ] (1,0) circle(1);
\draw[fill=blue!10,opacity=0.2 ] (3,0) circle(1);
\draw[fill=blue!10,opacity=0.2 ] (5,0) circle(1);
\draw[fill=blue!10,opacity=0.2 ] (7,0) circle(1);
\draw[fill=blue!10,opacity=0.2 ] (9,0) circle(1);
\draw[fill=blue!10,opacity=0.2 ] (2,1.75) circle(1);
\end{tikzpicture}
\quad
\begin{tikzpicture}[scale=0.4, line cap=round]
\draw[dashed] (0,0)--(10,0);
\renewcommand{\accX}{0}
\renewcommand{\accY}{0}
\po\NE\po\E\po\EB
\po\E\po\EB\po
\NEB\po\NE
\po
\E\po
\EB\po\E\po
\SE\po
\SEB
\po\SE\po
\draw[fill=blue!10,opacity=0.2 ] (1,0) circle(1);
\draw[fill=blue!10,opacity=0.2 ] (3,0) circle(1);
\draw[fill=blue!10,opacity=0.2 ] (5,0) circle(1);
\draw[fill=blue!10,opacity=0.2 ] (9,0) circle(1);
\draw[fill=blue!10,opacity=0.2 ] (6,1.75) circle(1);
\draw[fill=blue!10,opacity=0.2 ] (8,1.75) circle(1);
\end{tikzpicture}
\quad
\begin{tikzpicture}[scale=0.4, line cap=round]
\draw[dashed] (0,0)--(10,0);
\renewcommand{\accX}{0}
\renewcommand{\accY}{0}
\po\NE\po\E\po\EB
\po\NEB\po\NE\po\E\po\EB\po
\E\po\SE\po\SEB
\po
\EB\po\E

\po\SE\po
\draw[fill=blue!10,opacity=0.2 ] (1,0) circle(1);
\draw[fill=blue!10,opacity=0.2 ] (3,0) circle(1);
\draw[fill=blue!10,opacity=0.2 ] (7,0) circle(1);
\draw[fill=blue!10,opacity=0.2 ] (9,0) circle(1);
\draw[fill=blue!10,opacity=0.2 ] (6,1.75) circle(1);
\draw[fill=blue!10,opacity=0.2 ] (4,1.75) circle(1);
\end{tikzpicture}
\quad
\begin{tikzpicture}[scale=0.4, line cap=round]
\draw[dashed] (0,0)--(10,0);
\renewcommand{\accX}{0}
\renewcommand{\accY}{0}
\po\NE\po\NEB\po\NE
\po\E\po\EB\po\E\po\SE\po\SEB\po
\EB\po\E\po\EB\po\E\po\SE
\po
\draw[fill=blue!10,opacity=0.2 ] (1,0) circle(1);
\draw[fill=blue!10,opacity=0.2 ] (2,1.75) circle(1);
\draw[fill=blue!10,opacity=0.2 ] (7,0) circle(1);
\draw[fill=blue!10,opacity=0.2 ] (9,0) circle(1);
\draw[fill=blue!10,opacity=0.2 ] (5,0) circle(1);
\draw[fill=blue!10,opacity=0.2 ] (4,1.75) circle(1);
\end{tikzpicture}
\quad  
\begin{tikzpicture}[scale=0.4, line cap=round]
\draw[dashed] (0,0)--(10,0);
\renewcommand{\accX}{0}
\renewcommand{\accY}{0}
\po\NE\po\NEB\po\NE
\po\E\po\EB\po\E\po\EB\po\E
\SE\po\SEB\po
\EB\po\E\po\SE
\po
\draw[fill=blue!10,opacity=0.2 ] (1,0) circle(1);
\draw[fill=blue!10,opacity=0.2 ] (2,1.75) circle(1);
\draw[fill=blue!10,opacity=0.2 ] (7,0) circle(1);
\draw[fill=blue!10,opacity=0.2 ] (9,0) circle(1);
\draw[fill=blue!10,opacity=0.2 ] (6,1.75) circle(1);
\draw[fill=blue!10,opacity=0.2 ] (4,1.75) circle(1);
\end{tikzpicture}
\begin{tikzpicture}[scale=0.4, line cap=round]
\draw[dashed] (0,0)--(10,0);
\renewcommand{\accX}{0}
\renewcommand{\accY}{0}
\po\NE\E\po\EB\po\NEB\po\NE
\po\E\po\EB\po\E\po\EB\po\E
\SE\po\SEB\po
\SE
\po
\draw[fill=blue!10,opacity=0.2 ] (1,0) circle(1);
\draw[fill=blue!10,opacity=0.2 ] (3,0) circle(1);
\draw[fill=blue!10,opacity=0.2 ] (8,1.75) circle(1);
\draw[fill=blue!10,opacity=0.2 ] (9,0) circle(1);
\draw[fill=blue!10,opacity=0.2 ] (6,1.75) circle(1);
\draw[fill=blue!10,opacity=0.2 ] (4,1.75) circle(1);
\end{tikzpicture}
\quad 
\begin{tikzpicture}[scale=0.4, line cap=round]
\draw[dashed] (0,0)--(10,0);
\renewcommand{\accX}{0}
\renewcommand{\accY}{0}
\po\NE\po\NEB\po\NE\po\E\po\EB
\po\E\po\EB\po\E\po\EB\po\E
\SE\po\SEB\po
\SE
\po
\draw[fill=blue!10,opacity=0.2 ] (1,0) circle(1);
\draw[fill=blue!10,opacity=0.2 ] (2,1.75) circle(1);
\draw[fill=blue!10,opacity=0.2 ] (8,1.75) circle(1);
\draw[fill=blue!10,opacity=0.2 ] (9,0) circle(1);
\draw[fill=blue!10,opacity=0.2 ] (6,1.75) circle(1);
\draw[fill=blue!10,opacity=0.2 ] (4,1.75) circle(1);
\end{tikzpicture}
\quad 
\begin{tikzpicture}[scale=0.4, line cap=round]
\draw[dashed] (0,0)--(8,0);
\renewcommand{\accX}{0}
\renewcommand{\accY}{0}
\po\NE\po\NEB\po\NE\po\NEB\po\NE\po\E\po
\SE\po\SEB\po\EB\po\E\po\SE
\SEB\po\SE\po
\draw[fill=blue!10,opacity=0.2 ] (1,0) circle(1);
\draw[fill=blue!10,opacity=0.2 ] (2,1.75) circle(1);
\draw[fill=blue!10,opacity=0.2 ] (6,1.75) circle(1);
\draw[fill=blue!10,opacity=0.2 ] (7,0) circle(1);
\draw[fill=blue!10,opacity=0.2 ] (4,1.75) circle(1);
\draw[fill=blue!10,opacity=0.2 ] (3,3.5) circle(1);
\end{tikzpicture}
\quad
\begin{tikzpicture}[scale=0.4, line cap=round]
\draw[dashed] (0,0)--(8,0);
\renewcommand{\accX}{0}
\renewcommand{\accY}{0}
\po\NE\po\NEB\po\NE\po\E\po\EB\po\NEB\po\NE\po\E\po
\SE\po\SEB\po\SE
\SEB\po\SE\po
\draw[fill=blue!10,opacity=0.2 ] (1,0) circle(1);
\draw[fill=blue!10,opacity=0.2 ] (2,1.75) circle(1);
\draw[fill=blue!10,opacity=0.2 ] (6,1.75) circle(1);
\draw[fill=blue!10,opacity=0.2 ] (7,0) circle(1);
\draw[fill=blue!10,opacity=0.2 ] (4,1.75) circle(1);
\draw[fill=blue!10,opacity=0.2 ] (5,3.5) circle(1);
\end{tikzpicture}
\quad
\begin{tikzpicture}[scale=0.4, line cap=round]
\draw[dashed] (0,0)--(8,0);
\renewcommand{\accX}{0}
\renewcommand{\accY}{0}
\po\NE\po\NEB\po\NE\po\NEB\po\NE\po\E\po\EB\po\E\po
\SE\po\SEB\po\SE
\SEB\po\SE\po
\draw[fill=blue!10,opacity=0.2 ] (1,0) circle(1);
\draw[fill=blue!10,opacity=0.2 ] (2,1.75) circle(1);
\draw[fill=blue!10,opacity=0.2 ] (6,1.75) circle(1);
\draw[fill=blue!10,opacity=0.2 ] (7,0) circle(1);
\draw[fill=blue!10,opacity=0.2 ] (3,3.5) circle(1);
\draw[fill=blue!10,opacity=0.2 ] (5,3.5) circle(1);
\end{tikzpicture}
\quad
\begin{tikzpicture}[scale=0.4, line cap=round]
\draw[dashed] (0,0)--(8,0);
\renewcommand{\accX}{0}
\renewcommand{\accY}{0}
\po\NE\po\NEB\po\NE\po\E\po\SE\po\SEB\po\EB\po

\po\NEB\po\NE\po\E
\po\SE\po\SEB\po\SE\po

\draw[fill=blue!10,opacity=0.2 ] (1,0) circle(1);
\draw[fill=blue!10,opacity=0.2 ] (3,0) circle(1);
\draw[fill=blue!10,opacity=0.2 ] (5,0) circle(1);
\draw[fill=blue!10,opacity=0.2 ] (7,0) circle(1);
\draw[fill=blue!10,opacity=0.2 ] (6,1.75) circle(1);
\draw[fill=blue!10,opacity=0.2 ] (2,1.75) circle(1);
\end{tikzpicture}
\quad
\begin{tikzpicture}[scale=0.4, line cap=round]
\draw[dashed] (0,0)--(8,0);
\renewcommand{\accX}{0}
\renewcommand{\accY}{0}
\po\NE\po\NEB\po\NE\po\NEB\po\NE\po\E\po
\SE\po\SEB\po\SE\po\SEB\po
\EB\po\E\po\SE\po
\draw[fill=blue!10,opacity=0.2 ] (1,0) circle(1);
\draw[fill=blue!10,opacity=0.2 ] (2,1.75) circle(1);
\draw[fill=blue!10,opacity=0.2 ] (5,0) circle(1);
\draw[fill=blue!10,opacity=0.2 ] (7,0) circle(1);
\draw[fill=blue!10,opacity=0.2 ] (4,1.75) circle(1);
\draw[fill=blue!10,opacity=0.2 ] (3,3.5) circle(1);
\end{tikzpicture}
\quad
\begin{tikzpicture}[scale=0.4, line cap=round]
\draw[dashed] (0,0)--(8,0);
\renewcommand{\accX}{0}
\renewcommand{\accY}{0}
\po\NE\po\E\po\EB\po\NEB\po\NE\po\NEB\po\NE\po\E\po
\SE\po\SEB\po\SE
\SEB\po\SE\po
\draw[fill=blue!10,opacity=0.2 ] (1,0) circle(1);
\draw[fill=blue!10,opacity=0.2 ] (3,0) circle(1);
\draw[fill=blue!10,opacity=0.2 ] (6,1.75) circle(1);
\draw[fill=blue!10,opacity=0.2 ] (7,0) circle(1);
\draw[fill=blue!10,opacity=0.2 ] (4,1.75) circle(1);
\draw[fill=blue!10,opacity=0.2 ] (5,3.5) circle(1);
\end{tikzpicture}
\quad
\caption{The 17 inchworm paths with a kissing number equal to 6.}
    \label{fig:inchwormpath}
\end{figure}
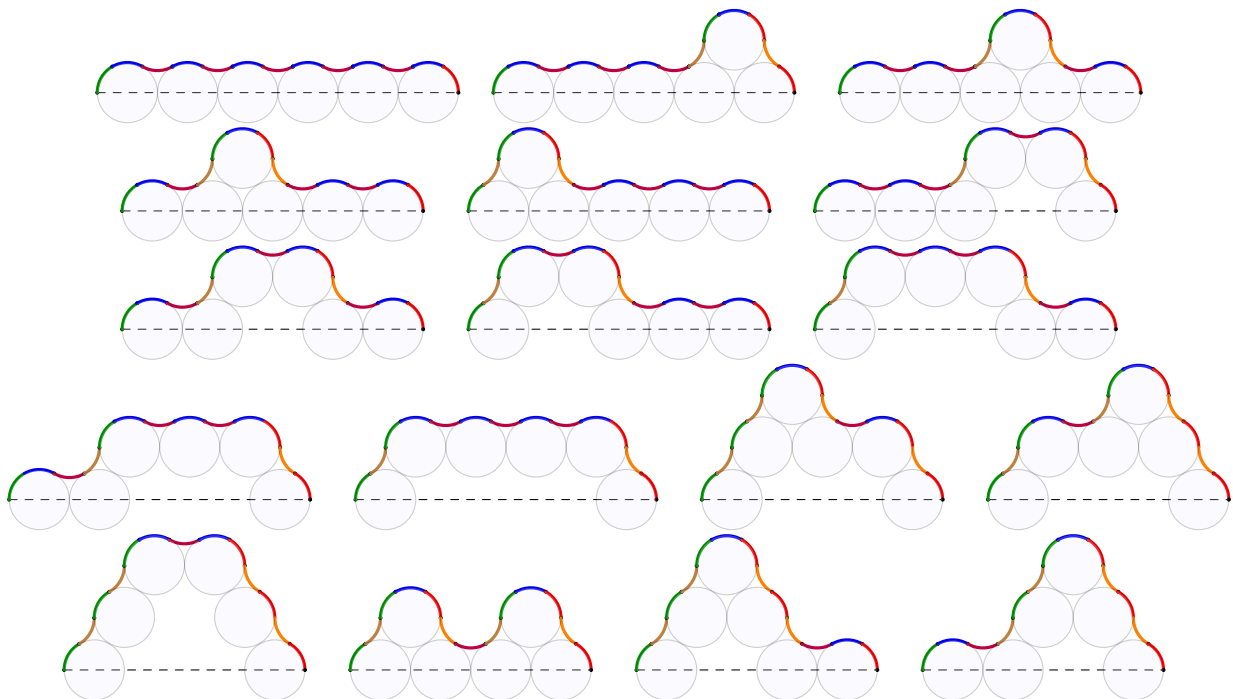

Now we consider the packing paths in $\mathcal{A}$ that avoid the two patterns $\overline{U}\overline{D}$ and $DU$. 
Note that such paths, with a given kissing number $k$, are in one-to-one correspondence with sequences $C=C_1C_2\cdots C_k$, where each  $C_i$ ($1\leq i\leq k$) represents a circle of radius two satisfying the following conditions
\begin{enumerate}
    \item  $C_1$ is centered at $(2,0)$;
    \item  $C_k$ is centered on the $x$-axis;
    \item  for every $i\geq 2$, the circle $C_i$ is tangent to $C_{i-1}$ according to the hexagonal lattice, and the center of $C_i$ is always on the right of $C_{i-1}$.
\end{enumerate}
Because of their resemblance to the movement of a caterpillar, these paths will be referred to as \emph{inchworm paths}. See Figure~\ref{fig:inchwormpath} for an illustration of the 17 inchworm paths with a kissing number equal to 6.

Let $A'(x, y)$ (resp.\ $B'(x, y)$) be the bivariate generating function that counts the number of inchworm paths in $\mathcal{A}$ (resp.\ $\mathcal{B}$), where  $x$ marks the width and $y$ marks the kissing number. 

Using the same decomposition with respect to its first return on the $x$ axis, we obtain the two functional equations 
$$A'(x,y) = 1+x^2B'(x,y), \quad  \mbox{  and }$$
$$B'(x,y) =  x^2y+x^4yB'(x,y)+
  x^2y^2(A'(x,y)-1)+x^4y^2B'(x,y)(A'(x,y)-1).$$

Consequently, we have the following.
\begin{thm} The generating function for the number of inchworm paths with respect to the width and the kissing number is given by \[
A'(x,y)
=\frac{1 - x^4 y + x^4 y^2 -\sqrt{1 - 2 x^4 y - 2 x^4 y^2 + x^8 y^2 - 2 x^8 y^3 + x^8 y^4}}{2 x^4 y^2}.
\]
\end{thm}
The first few terms of the power series expansion are
\begin{align*} 
1 + yx^4 &+ (y^3 + y^2)x^8 + (y^5 + 3y^4 + y^3)x^{12} + (y^7 + \bm{6y^6} + 6y^5 + y^4)x^{16} +\\
&(y^9 + 10y^8 + 20y^7 + \bm{10y^6} + y^5)x^{20} + (y^{11} + 15y^{10} + 50y^9 + 50y^8 + 15y^7 + \bm{y^6})x^{24} + O(x^{28}).
\end{align*}

Figure \ref{fig:inchwormpath} shows the inchworm paths with kissing number equal to 6, corresponding to the three bold
coefficients in the above series.

Consider the matrix $\mathcal{M}=[a(4n,n+k-1)]_{n, k\geq 1}$,  where $a(n,k)$ is the coefficient of $x^ny^k$ in the series expansion of $A'(x,y)$.  The first few rows of $\mathcal{M}$ are
\[
\footnotesize 
\mathcal{M}=\left(
\begin{array}{ccccccc}
 1 & 0 & 0 & 0 & 0 & 0 & 0 \\
 1 & 1 & 0 & 0 & 0 & 0 & 0 \\
 1 & 3 & 1 & 0 & 0 & 0 & 0 \\
 1 & 6 & 6 & 1 & 0 & 0 & 0 \\
 1 & 10 & 20 & 10 & 1 & 0 & 0 \\
 1 & 15 & 50 & 50 & 15 & 1& 0 \\
\end{array} 
\right), 
\]
which coincides with the Narayana triangle \oeis{A001263}. Therefore we have
$$a(4n,n+k-1)= \frac{1}{k} {n-1\choose k-1}{n\choose k-1}. \quad n, k \geq 1.$$

Observe that the number of inchworm paths of width $4n$ is the Catalan number   $c_{n}=\frac{1}{n+1}{2n\choose n}$. Let us define a bijection $f$ from inchworm paths of width $4n$ to Dyck paths of length $2n$ that transports the kissing number $k$ into $2n-k$ peaks $UD$.

An inchworm path can be characterized by a sequence of $C=C_1C_2\cdots C_k$ where $C_i$, $1\leq i\leq k$,  are circles satisfying three conditions introduced at the beginning of this section. So, we can decompose $C$ into two forms: ($a$) $C=C_1 C'$ where $C'$ is a sequence of circles corresponding to an inchworm path; and ($b$) $C=C_1 C' C_k C''$ where $C'$ and $C''$  are two sequences of circles corresponding to inchworm paths such that, $C'$ is the translation of a nonempty inchworm path by the vector $(2,2\sqrt{3})$, and  $C''$ starts and ends on the $x$-axis. From this decomposition, we can naturally  construct the  bijection $f$ as follows:
\[
f(\varepsilon) = \varepsilon, \quad  
f(C) =
\begin{cases}
UD f(C'), & \text{if } C = C_1 C' \text{ (form ($a$))}, \\[4pt]
U f(C') D f(C''), & \text{if } C = C_1 C' C_k C'' \text{ (form ($b$))}.
\end{cases}
\]

We can check that:

- If $C$ is of form ($a$) then $$\nbp(f(C))=1+\nbp(f(C'))=1+2(n-1)-\ki(C')=2n-\ki(C).$$

-  If $C$ is of form ($b$) then $$\nbp(f(C))=\nbp(f(C'))+\nbp(f(C''))=2n_1-\ki(C')+2n_2-\ki(C''),$$
 where  $4n_1=\s(C')$, $4n_2=\s(C'')$ and $n_1+n_2+1=n.$
So, we have \[\nbp(f(C))=2n-2-\ki(C')-\ki(C'')=2n-\ki(C).\]
\begin{cor} The generating function for the number of inchworm paths with respect to the kissing number is given by 
$$A'(1,y)=\frac{1 - y + y^2 -\sqrt{1 - 2 y - y^2 - 2 y^3 + y^4}}{2y^2}.$$
\end{cor}

The first terms of the series expansion of $A(1,y)$ are
$$1+y+y^2+2y^3+4y^4+8y^5+17y^6+37y^7+82y^8+185y^9+423y^{10}+O(y^{11}).$$
The sequence of coefficients corresponds to twice the sequence \oeis{A004148}. This sequence counts the number of secondary structures of RNA molecules, or equivalently the number of peakless Motzkin paths (Motzkin paths avoiding peak $UD$). Indeed, an inchworm path can be characterized by a sequence of $C=C_1C_2\cdots C_k$ as previously described. So, we can decompose $C$ into the  two same forms ($a$) and ($b$)  used for the bijection $f$. From this decomposition, we can naturally  construct a bijection $g$ from these paths to peakless Motzkin paths as follows:
\[
g(\varepsilon) = \varepsilon, \quad  
g(C) =
\begin{cases}
F g(C'), & \text{if } C = C_1 C' \text{ (form (a))}, \\[4pt]
U g(C') D g(C''), & \text{if } C = C_1 C' C_k C'' \text{ (form (b))}.
\end{cases}
\]

 \section{Going further}
Paths have been extensively investigated in the literature, both in the context of pattern avoidance and in the enumeration according to the number of pattern occurrences. In a future study, we plan to  carry out an exhaustive study of the distribution and popularity of patterns of size at most three. Such an analysis would likely lead to the discovery of new bijections between these paths and more classical combinatorial objects.

\begin{figure}[ht!]
    \centering
\newcommand{\ChainArcRev}[3]{%
  \CalcArcPoints{#1}{#2}%
  \MinMaxX{\startX}{\startY}{\endX}{\endY}%
  \pgfmathsetmacro{\dx}{\accX - \maxX}
  \pgfmathsetmacro{\dy}{\accY - \maxY}
  \begin{scope}[shift={(\dx,\dy)}]
    \draw[#3, very thick] ({cos(#1)},{sin(#1)}) arc[start angle=#1, end angle=#2, radius=\R];
  \end{scope}
  \pgfmathsetmacro{\accX}{\minX + \dx}
  \pgfmathsetmacro{\accY}{\minY + \dy}
}
\def\NErev{\ChainArcRev{120}{180}{green!60!black}}  
\def\NEBrev{\ChainArcRev{300}{360}{brown}}          

    \begin{tikzpicture}[scale=0.4, line cap=round]
      \draw[dashed] (0,0)--(16,0);
      \renewcommand{\accX}{0}
      \renewcommand{\accY}{0}
      
      \filldraw[black] (0,0) circle (0.06);
      
      \NE   \filldraw[black] (\accX,\accY) circle (0.05);
      \NEB  \filldraw[black] (\accX,\accY) circle (0.05);
            \NE   \filldraw[black] (\accX,\accY) circle (0.05);
             \E    \filldraw[black] (\accX,\accY) circle (0.05);
             \SE   \filldraw[black] (\accX,\accY) circle (0.05);
               \NEBrev  \filldraw[black] (\accX,\accY) circle (0.05);
         \E   \filldraw[black] (\accX,\accY) circle (0.05);
         \EB  \filldraw[black] (\accX,\accY) circle (0.05);
        \NEB   \filldraw[black] (\accX,\accY) circle (0.05);
            \SEB  \filldraw[black] (\accX,\accY) circle (0.05);
             \EB  \filldraw[black] (\accX,\accY) circle (0.05);
    \NEB  \filldraw[black] (\accX,\accY) circle (0.05);
      \NE   \filldraw[black] (\accX,\accY) circle (0.05);
    \E   \filldraw[black] (\accX,\accY) circle (0.05);
        \EB   \filldraw[black] (\accX,\accY) circle (0.05);
              \NEB  \filldraw[black] (\accX,\accY) circle (0.05);
        \NE   \filldraw[black] (\accX,\accY) circle (0.05);
            \E   \filldraw[black] (\accX,\accY) circle (0.05);
          \SE  \filldraw[black] (\accX,\accY) circle (0.05);
\NEBrev \filldraw[black] (\accX,\accY) circle (0.05);
\NErev \filldraw[black] (\accX,\accY) circle (0.05);
         \SEB  \filldraw[black] (\accX,\accY) circle (0.05);  
                \EB   \filldraw[black] (\accX,\accY) circle (0.05);
                   \E   \filldraw[black] (\accX,\accY) circle (0.05);
                      \EB   \filldraw[black] (\accX,\accY) circle (0.05);
                       \NErev   \filldraw[black] (\accX,\accY) circle (0.05);
      \draw[fill=blue!10,opacity=0.2 ] (1,0) circle(1);
      \draw[fill=blue!10,opacity=0.2 ] (2,1.75) circle(1);
      \draw[fill=blue!10,opacity=0.2 ] (3,0) circle(1);
          \draw[fill=blue!10,opacity=0.2 ] (4,1.75) circle(1);
      \draw[fill=blue!10,opacity=0.2 ] (5,0) circle(1);
          \draw[fill=blue!10,opacity=0.2 ] (6,1.75) circle(1);
      \draw[fill=blue!10,opacity=0.2 ] (7,0) circle(1);
          \draw[fill=blue!10,opacity=0.2 ] (8,1.75) circle(1);
      \draw[fill=blue!10,opacity=0.2 ] (9,0) circle(1);
      \draw[fill=blue!10,opacity=0.2 ] (10,1.75) circle(1);
      \draw[fill=blue!10,opacity=0.2 ] (11,0) circle(1);
      \draw[fill=blue!10,opacity=0.2 ] (12,1.75) circle(1);
      \draw[fill=blue!10,opacity=0.2 ] (13,0) circle(1);
      \draw[fill=blue!10,opacity=0.2 ] (14,1.75) circle(1);
      \draw[fill=blue!10,opacity=0.2 ] (15,0) circle(1);
       \draw[fill=blue!10,opacity=0.2 ] (9,3.5) circle(1);
        \draw[fill=blue!10,opacity=0.2 ] (11,3.5) circle(1);
         \draw[fill=blue!10,opacity=0.2 ] (13,3.5) circle(1);
    \end{tikzpicture}
    \caption{The skew packing  path $U\overline{U}UFD\bar{U}_LF\overline{F}\overline{U}\overline{D}\overline{F}\overline{U}UF\overline{F}\overline{U}UFD\bar{U}_LU_L\overline{D}\overline{F}F\overline{F}U_L$ of width 28, with 26  steps, of area~11  and kissing number~10.}
    \label{fig:pathskew}
\end{figure}
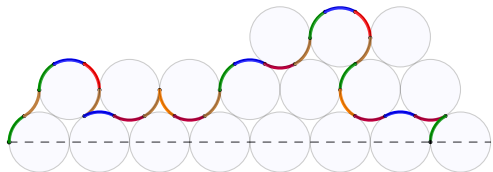

A \emph{skew packing path} is a self-avoiding path in the hexagonal circle packing starting at the origin, ending on the $x$-axis, and consisting of steps from $S=\{U,F,D, \bar{U}, \bar{F}, \bar{D}\}$, and where the two steps $U$ and $\bar{U}$ can be crossed from right to left or  from left to right (while maintaining self-avoidance). We denote by $U_L$ and $\bar{U}_L$ these steps crossed from right to left. See Figure~\ref{fig:pathskew} for an example of such a path. 
 Can one obtain multivariate generating functions for the enumeration of skew packing paths with respect to the parameters studied in this paper (number of steps, width, area, kissing number)?

\begin{figure}[H]
    \centering
    \begin{tikzpicture}[scale=0.4, line cap=round]
      \draw[dashed] (0,0)--(16,0);
      \renewcommand{\accX}{0}
      \renewcommand{\accY}{0}
      
      \filldraw[black] (0,0) circle (0.06);
      
      \NE   \filldraw[black] (\accX,\accY) circle (0.05);
      \NEB  \filldraw[black] (\accX,\accY) circle (0.05);
      \SEB  \filldraw[black] (\accX,\accY) circle (0.05);
      \EB   \filldraw[black] (\accX,\accY) circle (0.05);
      \E    \filldraw[black] (\accX,\accY) circle (0.05);
      \SE   \filldraw[black] (\accX,\accY) circle (0.05);
      \NE   \filldraw[black] (\accX,\accY) circle (0.05);
      \E    \filldraw[black] (\accX,\accY) circle (0.05);
      \EB   \filldraw[black] (\accX,\accY) circle (0.05);
      \NEB  \filldraw[black] (\accX,\accY) circle (0.05);
      \SEB  \filldraw[black] (\accX,\accY) circle (0.05);
      \EB   \filldraw[black] (\accX,\accY) circle (0.05);
      \NEB  \filldraw[black] (\accX,\accY) circle (0.05);
      \NE   \filldraw[black] (\accX,\accY) circle (0.05);
      \E    \filldraw[black] (\accX,\accY) circle (0.05);
       \EB   \filldraw[black] (\accX,\accY) circle (0.05);
       \NEB  \filldraw[black] (\accX,\accY) circle (0.05);
         \renewcommand{\accX}{5.5}
      \renewcommand{\accY}{2.6}
      
      \E    \filldraw[black] (\accX,\accY) circle (0.05);
       \EB    \filldraw[black] (\accX,\accY) circle (0.05);
        \NEB    \filldraw[black] (\accX,\accY) circle (0.05);
         \NE    \filldraw[black] (\accX,\accY) circle (0.05);
          \E    \filldraw[black] (\accX,\accY) circle (0.05);
           \SE    \filldraw[black] (\accX,\accY) circle (0.05);
            \NE    \filldraw[black] (\accX,\accY) circle (0.05);
             \E    \filldraw[black] (\accX,\accY) circle (0.05);
              \SE    \filldraw[black] (\accX,\accY) circle (0.05);

             \renewcommand{\accX}{5.5}
      \renewcommand{\accY}{2.6}
       \NEB    \filldraw[black] (\accX,\accY) circle (0.05);
        \NE    \filldraw[black] (\accX,\accY) circle (0.05);
        \NEB    \filldraw[black] (\accX,\accY) circle (0.05);
         \NE    \filldraw[black] (\accX,\accY) circle (0.05);
          \E    \filldraw[black] (\accX,\accY) circle (0.05);
           \EB    \filldraw[black] (\accX,\accY) circle (0.05);
            \E    \filldraw[black] (\accX,\accY) circle (0.05);
             \SE    \filldraw[black] (\accX,\accY) circle (0.05);
              \NE    \filldraw[black] (\accX,\accY) circle (0.05);
               \E    \filldraw[black] (\accX,\accY) circle (0.05);
             \SE    \filldraw[black] (\accX,\accY) circle (0.05);
             \SEB    \filldraw[black] (\accX,\accY) circle (0.05);
                \SE    \filldraw[black] (\accX,\accY) circle (0.05);
                 \SEB    \filldraw[black] (\accX,\accY) circle (0.05);
                \SE    \filldraw[black] (\accX,\accY) circle (0.05);
                   \renewcommand{\accX}{10.5}
      \renewcommand{\accY}{0.85}
       \E    \filldraw[black] (\accX,\accY) circle (0.05);
        \EB    \filldraw[black] (\accX,\accY) circle (0.05);
         \NEB    \filldraw[black] (\accX,\accY) circle (0.05);
         \SEB    \filldraw[black] (\accX,\accY) circle (0.05);
          \EB    \filldraw[black] (\accX,\accY) circle (0.05);
           \NEB    \filldraw[black] (\accX,\accY) circle (0.05);
      \draw[fill=blue!10,opacity=0.2 ] (1,0) circle(1);
      \draw[fill=blue!10,opacity=0.2 ] (2,1.75) circle(1);
      \draw[fill=blue!10,opacity=0.2 ] (3,0) circle(1);
          \draw[fill=blue!10,opacity=0.2 ] (4,1.75) circle(1);
      \draw[fill=blue!10,opacity=0.2 ] (5,0) circle(1);
          \draw[fill=blue!10,opacity=0.2 ] (6,1.75) circle(1);
      \draw[fill=blue!10,opacity=0.2 ] (7,0) circle(1);
          \draw[fill=blue!10,opacity=0.2 ] (8,1.75) circle(1);
      \draw[fill=blue!10,opacity=0.2 ] (9,0) circle(1);
      \draw[fill=blue!10,opacity=0.2 ] (10,1.75) circle(1);
      \draw[fill=blue!10,opacity=0.2 ] (11,0) circle(1);
      \draw[fill=blue!10,opacity=0.2 ] (12,1.75) circle(1);
      \draw[fill=blue!10,opacity=0.2 ] (13,0) circle(1);
      \draw[fill=blue!10,opacity=0.2 ] (14,1.75) circle(1);
      \draw[fill=blue!10,opacity=0.2 ] (15,0) circle(1);
       \draw[fill=blue!10,opacity=0.2 ] (13,3.5) circle(1);
        \draw[fill=blue!10,opacity=0.2 ] (5,3.5) circle(1);
         \draw[fill=blue!10,opacity=0.2 ] (7,3.5) circle(1);
          \draw[fill=blue!10,opacity=0.2 ] (9,3.5) circle(1);
           \draw[fill=blue!10,opacity=0.2 ] (11,3.5) circle(1);
            \draw[fill=blue!10,opacity=0.2 ] (8,5.25) circle(1);
            \draw[fill=blue!10,opacity=0.2 ] (10,5.25) circle(1);      \draw[fill=blue!10,opacity=0.2 ] (12,5.25) circle(1);
    \end{tikzpicture}
    \caption{A self-avoiding packing path in the quarter plane.}
    \label{fig:selfavoid}
\end{figure}
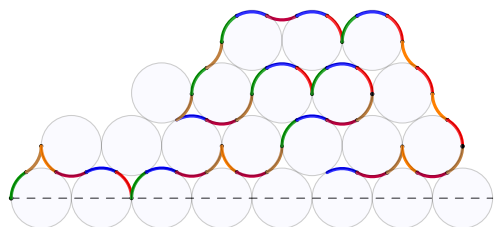

Finally, enumerating self-avoiding paths is a long-standing and challenging problem in combinatorics. A deeper understanding of their structure on this circle packing may uncover new bijective relationships and extend classical enumeration results for these   paths. See Figure~\ref{fig:selfavoid} for an example of such a path.

\end{document}